\documentclass[a4paper,10pt,reqno,final]{amsart}

\usepackage{epsfig}
\usepackage{psfrag}

\usepackage{amsmath}
\usepackage{amsfonts}
\usepackage{amsthm}

\usepackage[T1]{fontenc}

\usepackage[latin1]{inputenc}

\usepackage[french,english]{babel}

\usepackage[centerlast,small]{subfigure}

\usepackage{hyperref}
\hypersetup{
breaklinks=true,
}

\theoremstyle{definition}

\theoremstyle{remark}

\numberwithin{equation}{section}

\newcommand{\Er}{\mathbb{R}}

\newcommand{\vnorm}[2][{}]{\left\Vert#2\right\Vert_{#1}}

\makeatletter
\renewcommand{\@seccntformat}[1]{\large{\csname the#1\endcsname}.
\hspace{0.5em}}
\renewcommand{\section}{\@startsection {section}{1}{0mm}%
                                   {-\baselineskip}%
                                   {0.5\baselineskip}%
                                   {\sffamily\large\upshape\bfseries}}
\renewcommand{\subsection}{\@startsection {subsection}{2}{0mm}%
                                   {-0.5\baselineskip}%
                                   {0.5\baselineskip}%
                                   {\sffamily\normalsize\upshape
                                   \bfseries}}
\makeatother

\begin{document}
\selectlanguage{english}



\title%
[Estimation of anthropometrical and inertial body parameters]%
{Estimation of anthropometrical and inertial body parameters using double integration of residual torques and forces during squat jump}

\author{J\'er\^ome Bastien}
\author{Yoann Blache}
\author{Karine Monteil}

\date{\today}

\selectlanguage{french}

\address{
Université de Lyon\\
Centre de Recherche et d'Innovation sur le Sport\\
   U.F.R.S.T.A.P.S.\\
   Université Claude Bernard - Lyon 1\\
   27-29, Bd du 11 Novembre 1918\\
   69622 Villeurbanne Cedex\\
France}

\email{jerome.bastien@univ-lyon1.fr}
\email{yoann.blache@univ-lyon1.fr}
\email{karine.monteil@univ-lyon1.fr}

\selectlanguage{english}

\begin{abstract}
The inertial (IP) and anthropometrical (AP) parameters of human body are mostly estimated from coefficients issue from cadaver measurements. These parameters could involve errors in the calculation of joint torques during explosive movements. The purpose of this study  was to optimize the IP and AP in order to minimize the residual torque and force during squat jumping.  Three  methods
 of determination have been presented: 
method A:  optimizing AP and IP of each body part,
method B:  optimizing trunk AP and IP, assuming that the AP and IP of the lower limbs were known,
method C:  using Winter AP and IP.
For each method, the value (degree 0),  the integral (degree 1) and the double integral (degree 2) of the residual moment
were also used. The method B with degree 2 was the most accurate to determine trunk AP and IP by minimizing  of the
residual force and torque, by providing a linear least squares system. Instead of minimizing the residual force and torque, by
classical way, the double integral of the latter provided more accurate results.
\end{abstract}


\maketitle

\markboth{Jérôme BASTIEN et al.}
{\MakeUppercase{%
Estimation of anthropometrical and inertial body parameters}}


\newcommand{\alphaquatre}{0.4796}

\newcommand{\alphaunwinter}{0.5000}
\newcommand{\alphadeuxwinter}{0.5670}
\newcommand{\alphatroiswinter}{0.5670}
\newcommand{\alphaquatrewinter}{0.6260}
\newcommand{\mparunwinter}{0.0290}
\newcommand{\mpardeuxwinter}{0.0930}
\newcommand{\mpartroiswinter}{0.2000}
\newcommand{\mparquatrewinter}{0.6780}

\newcommand{\residuvalD}{0.49027758}
\newcommand{\coeffdetvalD}{-0.93202331}

\newcommand{\residuintD}{0.46254100}
\newcommand{\coeffdetintD}{-2.04019651}


\newcommand{\normeldeuxintmoresD}{54.95696888}

\newcommand{\residuintdD}{0.38002630}
\newcommand{\coeffdetintdD}{-0.95403457}

\newcommand{\residuvalC}{0.45798897}
\newcommand{\coeffdetvalC}{0.10704717}

\newcommand{\residuintC}{0.25915191}
\newcommand{\coeffdetintC}{0.56793975}


\newcommand{\normeldeuxintmoresC}{12.59295192}

\newcommand{\residuintdC}{0.12324873}
\newcommand{\coeffdetintdC}{0.89740132}

\newcommand{\residuvalA}{0.31124363}
\newcommand{\coeffdetvalA}{0.54983769}


\newcommand{\inertieunval}{-21.10836389}
\newcommand{\inertiedeuxval}{-9.67961630}
\newcommand{\inertietroisval}{-1.56382741}
\newcommand{\inertiequatreval}{-3.63374661}

\newcommand{\residuintA}{0.06599102}
\newcommand{\coeffdetintA}{0.98053164}


\newcommand{\inertieunint}{9.96412006}
\newcommand{\inertiedeuxint}{10.13963418}
\newcommand{\inertietroisint}{0.58549936}
\newcommand{\inertiequatreint}{8.28949873}


\newcommand{\normeldeuxintmoresA}{0.27828425}

\newcommand{\residuintdA}{0.00308765}
\newcommand{\coeffdetintdA}{0.99994990}


\newcommand{\inertieunintd}{13.13320316}
\newcommand{\inertiedeuxintd}{12.55758823}
\newcommand{\inertietroisintd}{1.29624104}
\newcommand{\inertiequatreintd}{9.52771715}

\newcommand{\residuvalB}{0.38646452}
\newcommand{\coeffdetvalB}{0.24688232}

\newcommand{\residuintB}{0.25223184}
\newcommand{\coeffdetintB}{0.69054004}


\newcommand{\normeldeuxintmoresB}{4.27317154}

\newcommand{\residuintdB}{0.04714447}
\newcommand{\coeffdetintdB}{0.98865321}

\newcommand{\rnormunwinter}{0.4750}
\newcommand{\rnormdeuxwinter}{0.3020}
\newcommand{\rnormtroiswinter}{0.3230}
\newcommand{\rnormquatrewinter}{0.4960}

\newcommand{\inertieunwinter}{0.0116}
\newcommand{\inertiedeuxwinter}{0.1104}
\newcommand{\inertietroiswinter}{0.2890}
\newcommand{\inertiequatrewinter}{4.2067}

\newcommand{\rnormquatre}{0.5587}


\section{Introduction}
\label{intro}

Joint forces and torques are commonly used in motion analysis for orthopedics, ergonomics or sports science \cite{Pearsall1994a,Reid1990}. 
A 
standard bottom-up inverse dynamic model is often used to calculate joint forces and torques of the lower extremity. The principle of bottom-up inverse dynamic is to combine kinematic data and ground reaction force to calculate    these parameters. Moreover, body anthropometric (AP) and inertial parameters (IP) are needed to apply inverse dynamic model. 
AP and IP can be obtained from many ways. Cadaver measurements have been the first method applied by researchers and is still commonly used 
\cite{Clauseretal1969,fujikawa1963,chandleretal1978,Hinrichs1990g,Dempster1955}. 
Then,  predictive linear or non-linear equation \cite{Hinrichs1985,zatsiorsk1983,McConville1980,winter,Yeadon1989} and imaging resonance magnetic techniques 
\cite{%
Durkin1998,
Huang1983,Martin1989,Pearsall1996,Cheng2000} have  been elaborated. In \cite{Hatze2002}, the authors  have been interested in the accuracy of the different methods. They observed that the methods using $\gamma$ ray, X ray, tomography and imaging resonance magnetic techniques presented an average accuracy of 5\% with a maximal error of 11\%. Linear regression yields AP and IP with an average accuracy of about 24\% with a maximal error of 40\%. The corresponding values for non-linear regression were 16\% and 38\% respectively. According to \cite{Hatze2002}, the most accurate technique would be the anthropometrico-computational 
method with 1.8\% of average accuracy and 3\% of maximal error.

Some authors tried to evaluate the influence of error in AP and IP on joint torques during gait analyses. In \cite{Goldberg2008}, the authors  compared joint torques during walking using cadaver AP and IP \cite{Dempster1955} versus direct measurements \cite{Fowler1999}. Inertial moment, mass segment and segment center of mass position were significantly different. During the stance phase, joint torques were similar, while significant joint torque differences were observed during the swing phase. In \cite{Pearsall1999},  6 methods to calculate AP and IP were compared. Even if these methods provided different AP and IP, no effect on torque measurements were pointed out during walking. Nevertheless, in \cite{Pearsall1999}, the authors concluded that changes in AP and IP should  have a greater influence on the torque measurements for activities involving greater accelerations.

Torque measurements, using a standard inverse dynamic routine, are also used in explosive movements such as vertical jumps. During squat jumping, the push-off lasts around 350 ms, and the acceleration of the body center mass could reach 20~$\text{m.}{\text{s}}^{-2}$. 
According to \cite{Pearsall1999,Plamondon1996}, the accuracy of AP and IP should be important in regard with the great acceleration found during vertical jumping. This movement is often studied in the 2D sagittal plane and inverse dynamic model is applied, most of the time, to calculate torque and work at the lower limb joints. Therefore AP and IP are needed. Researchers mainly use cadaver measurements \cite{Lees2000,Lees2004,Bobbert2006a,Bobbert2008a,Laffaye2005,Laffaye2007} or predictive equations \cite{Domire2007,Domire2010,Cheng2008,Hara2008,Wilson2007,Haguenauer2006,Vanrenterghem2004}.
However these methods give AP and IP which are not specific to the population studied. As a result, some errors in AP and IP could imply  inaccuracy in joint torque  calculations. 

Some authors \cite{Chen2011,Kingma1998}
 used static analysis based on the measurement of the center of pression in order to calculate 
center of mass.
{%
They measured the center of pression (via force plateform) during differents equilibrium standing positions which corresponded to the for-aft position of body 
the center of mass. Thus, considering the inertial and anthropemtric parameters of the lower limbs as known, the position the center of mass of the trunk segment 
was calculated.}
Other authors created methods to minimize the errors of 3-D inverse dynamic calculi. 
These methods consisted  to minimize the error between the ground reaction force measured and the ground reaction 
force calculated with a top-down inverse dynamic model \cite{Riemer2008,Riemer2009,Kuo1998,Vaughan1982}.
However, to the best of our knowledge, no study has focused on the optimization of AP and IP in explosive movements and especially in 2-D squat jumping investigation. Therefore the purpose of this study was to adjust AP and IP of the human segments  during squat jumping in order to minimize error in joint torque values. 

Especially, the optimization will focus on the "head arm trunk" segment (HAT). 
Unlike \cite{winter}, 
this segment is not considered as being rigid
 {as it is composed of three segments}. 
{%
Moreover,  the position of the arms, the head and the trunk may be differed  during squat jumping from the   position collected on cadavers.}

In the first method, IP of each body parts will be optimized.
The second  method consists in optimizing HAT IP only, assuming that the IP of the others body parts are  known.
The third  one refers to Winter's IP \cite{winter}.
Finally, for each method we tried to minimize the residual torque, the integral of residual torque and the double integral of the residual torque.
{This last method lead to the best result.}

\section{Methods}
\label{mat}

The following parts contain: in \ref{expe}, the  acquisition of the experimental data;
in \ref{mat_notation}, the  notations of the studied system;
in  \ref{syntrraireac}, the synchronization beetwen displacements and forces, and 
the determination of the  AP of the trunk (AP  of the others body parts being known);
in  \ref{mat_lissage}, the smoothing technique applied on the   experimental displacements;
in \ref{mat_DIF},
the  inverse dynamic procedures and the development of the three methods to estimate IP.

\subsection{Experimental acquisition}
\label{expe}

\begin{figure}[ht]
\psfrag{SJ}{Squat Jump}
\begin{center}
\epsfig{file=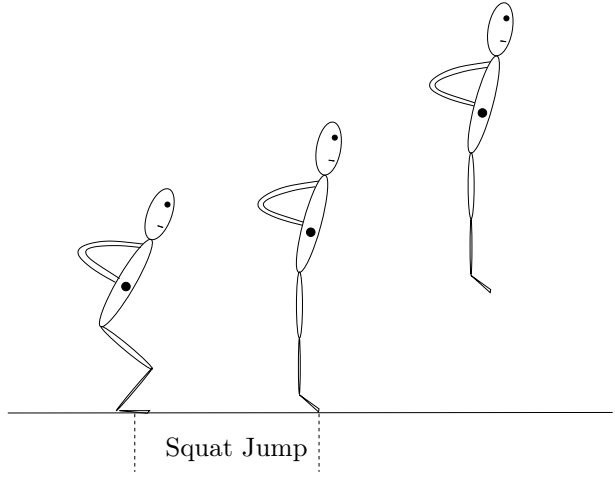,  width=8 cm}
\end{center}
\caption{\label{fig_SJ}Squat Jump.}
\end{figure}

Twelve healthy athletic male adults (mean $\pm$ SD: age, 23.2 $\pm$ 3.6 years; height, 1.75 $\pm$ 0.06 m; mass, 69.1 $\pm$ 8.2 kg) volunteered to participate in the study and provided informed consent. Prior to the experimental protocol, reflective landmarks were located on the right 5-th metatarsophalangeal, lateral malleolus, lateral femoral epicondyle, greater trochanter and acromion. A 10 minutes warm-up, including squat jumps session prepared participants to the task and allowed subjects to find their preferred squat depth position 
(see figure \ref{fig_SJ}).
Thereafter, the subjects performed at most ten maximal squat jumps. 
In order to avoid the contribution of the arms in vertical jump height \cite{Domire2010,Hara2008}, the subjects were instructed 
to keep their hands on their hip throughout the jump.
They also had to maintain the same initial squat depth for each jump. 

All jumps were performed on an AMTI force plate model OR6-7-2000 sampled at 1000 Hz. Countermovement defined as 
 a decrease of vertical ground reaction force ($R_y$) before the push-off phase was not allowed. 
 
       The beginning of the push-off was considered as the instant when the derivative of the smoothed $R_y$ is different to zero.
Simultaneously, the subjects were filmed in the sagittal plane with a 100 Hz camcorder (Ueye, IDS UI-2220SE-M-GL). The optical axis of the camcorder was perpendicular to the plane of the motion and located at 4 meters from the subject. 

      Jumps recorded were digitalized frame by frame with the Loco \copyright software (Paris, France). 
A four rigid segments model composed of the foot (left and right feet together), the shank (left and right shanks together), the thigh (left and right thighs together) and the HAT (head, arms and trunk) was used.
{%
Squat jump being a symmetrical motion, the lower limb segments were laterally combined  together and it was supposed that the left and right sides participate equivalently to the inter-articular efforts. The position of the upper limbs is fixed to limit their influence on $I_4$. Moreover the objective is to provide a robust estimate of the $I_4$ according to a given protocol and to observe that it is different from that of Winter.}

\subsection{Notations of the studied  system}
\label{mat_notation}

\begin{figure}[ht]
\psfrag{O}{$O$}
\psfrag{i}{$\vec i$}
\psfrag{j}{$\vec j$}
\psfrag{A1}{$A_1$}
\psfrag{A2}{$A_2$}
\psfrag{A3}{$A_3$}
\psfrag{Ap-2}{$A_{4}$}
\psfrag{Ap}{$A_5$}
\psfrag{t1}{$\theta_1$}
\psfrag{t2}{$\theta_2$}
\psfrag{t3}{$\theta_3$}
\psfrag{t4}{$\theta_4$}
\psfrag{X}{horizontal axis (X)}
\psfrag{Y}{vertical axis (Y)}
\begin{center}
\epsfig{file=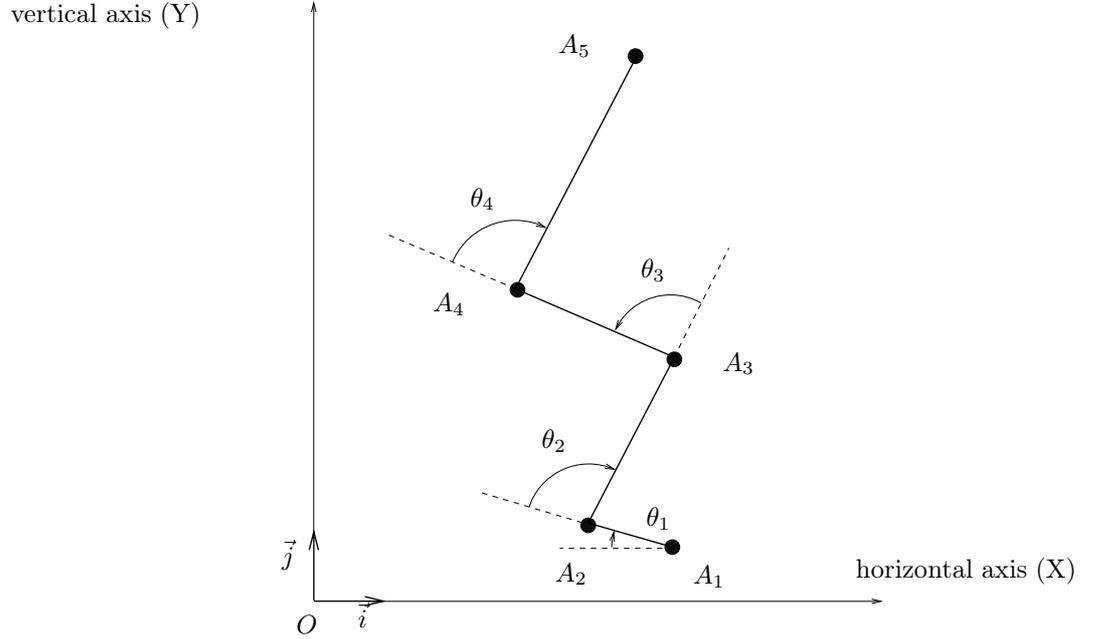,  width=11.5 cm}
\end{center}
\caption{\label{lec_kino_fig:29}Studied subject in his sagittal plan.
The considered anatomical landmarks are:
$A_1$, the 5-th metatarsophalangeal, 
$A_2$, the lateral malleolus, 
$A_3$, the lateral femoral epicondyle, 
$A_4$, the greater trochanter 
and $A_5$, the acromion. 
The joint angles are $\theta_1, \dots, \theta_4$.}
\end{figure}

For
\begin{equation}
\label{defq}
q=5,
\end{equation}
joint angles
$\theta_1, \dots, \theta_{q-1}$ are defined by 
\begin{subequations}
\label{eq01tot}
\begin{align}
\label{eq01a}
&\theta_1=\widehat{\left(\vec i, \overrightarrow{A_1A_{2}}\right)} \in (-\pi,\pi],\\
\label{eq01b}
&\forall j\in\{2,...,q-1\},\quad 
\theta_j=\widehat{\left(\overrightarrow{A_{j-1}A_{j}}, \overrightarrow{A_jA_{j+1}}\right)}\in (-\pi,\pi],
\end{align}
with $A_j$ the anatomical landmark, from the 5-th metatarsophalangeal to the acromion.
The constant lengths  $l_1, \dots, l_{q-1}$ are defined by 
\begin{equation}
\label{eq01c}
\forall j\in\{1,...,q-1\},\quad
l_j=A_jA_{j+1}. 
\end{equation}
\end{subequations}

The coordinates of points $A_j$, denoted $\left(x_j^i,y_j^i\right)$, 
was obtained from  experimental data: for all $j\in \{1,...,q\}$
\begin{equation}
\label{sjeq01}
\forall  i\in \{0,...,n\},\quad
x_j^i=x_j(i/f_e),\quad
y_j^i=y_j(i/f_e),
\end{equation}
and $f_e=100$ Hz was the acquisition frequency.
The center of mass of  segment $[A_j,A_{j+1}]$ is noted $G_j$ with  
\begin{equation*}
\alpha_j=\frac{A_jG_j}{A_jA_{j+1}},
\end{equation*}
$\alpha_j$
 being the distance from the distal joint to the segment center
of mass relative to the segment length: 
5-th metatarsophalangeal, 
lateral malleolus, 
lateral femoral epicondyle, 
the greater trochanter,
for the foot, the shank, the thigh an the trunk respectively
\cite{winter}.
$x_{G_j}^i$ and $y_{G_j}^i$ are the coordinates of $G_j$ at times $i/f_e$.
{As usually performed for kinematic and dynamic analysis of 2D movements in the sagittal plane 
\cite{Lees2004B,Bobbert2008a,Hara2008B},
the center 
of mass of each segment was assumed to lie on the line connecting the  markers.}
For each segment $\left[A_jA_{j+1}\right]$, anthropometry data are:
\begin{itemize}
\item[$\bullet$]
$m_j$, the mass;
\item[$\bullet$]
$l_j$, the length;
\item[$\bullet$]
$\alpha_j$, the relative position of center of mass $G_j$
\item[$\bullet$]
$I_j$, the moment of inertia according to its center of mass $G_j$.
\end{itemize}
The center of  mass of the subject is noted $G=(x_G,y_G)$ and his  total mass  $m$.
The data  $m_j$  and  
$\alpha_j$ for $1\leq j\leq 3$ and $m_4$ are  determined in \cite{winter}.
The value $\alpha_4$ will be determined by the method of Section \ref{syntrraireac}.
From length $l_j$ and mass $m_j$ of segment,
the value of the normalized radius of gyration $\widetilde  r_j$ is defined by 
$
\widetilde r_j=r_j/l_j,
$
where $r_j$ is the radius of gyration (according to the center of mass).
Then 
the moment of inertia is:
$
I_j=m_jr^2_j.
$
The values $I_j$ will be determined by the method of Section \ref{mat_DIF}.

For all, $1\leq j\leq q$, $\vec R_j=\left(R_{x,j},R_{y,j}\right)$ and $C_j$ are
respectively the resultant force and the torque, according to point $A_j$,  of the action of 
segment $\left[A_{j-1}A_{j}\right]$ on segment $\left[A_{j}A_{j+1}\right]$. Conventionally, for $j=1$, 
$\vec R_1=\vec R$ and  $C_1=C$ are the ground reaction (forces and torques, according to point $A_1$), and for 
$j=q$, $\vec R_q$ and $C_q$ are equal to zero.

The experimental data were the values  $\vec R$, and $C$, denoted $\vec R^i$ and $C^i$,  such as
\begin{equation}
\label{sjeq02}
\forall  i\in \{0,...,n'\},\quad
\vec R^i=\vec R(i/f'_e), \quad C^i=C(i/f'_e).
\end{equation}
and $f'_e=1000$ Hz was the acquisition frequency.
Actually displacement and force data are not synchronized, which means that 
 \eqref{sjeq02} 
has to be replaced by the following equation:
there exits an unknow integer $\nu$  such that 
\begin{equation}
\label{sjeqdecalage}
\vec R^i=\vec R\left(\frac{i+\nu}{f'_e}\right),\quad
C^i=\vec C\left(\frac{i+\nu}{f'_e}\right).
\end{equation}

\subsection{Synchronization of displacements and forces and determination of $\alpha_4$}
\label{syntrraireac}

The double integration of residual force is considered, this force is defined as the difference between
the measured experimental ground reaction force  and the theoretical  ground reaction force determined with respect to the position of the center of mass
{defined by \eqref{46tot}}.
The double integration depends on  the unknown integer  $\nu$ defined by \eqref{sjeqdecalage}
and the unknown  number $\alpha_4$, which corresponds to the position of the center of mass of the trunk.
By a double optimization on $\nu$ and $\alpha_4$, 
the norm  $l^2$ (also called Root Mean Square Error,  RMSE)
is minimal.
For more details, the reader should refer  to Appendix \ref{syntrraireacannexe}.

\subsection{Smoothing of experimental data $x_j^i$,  $y_j^i$, $x_{G_j}^i$, $y_{G_j}^i$, 
$x_{G}^i$ and $y_{G}^i$}
\label{mat_lissage}

Values of $x_j^i$,  $y_j^i$, $x_{G_j}^i$, $y_{G_j}^i$, 
$x_{G}^i$ and $y_{G}^i$, are experimental 
data.
They need to be smoothed in order to be 
 derivatived once or twice.
 The smoothing has for objective to minimize the residual reaction force. Nevertheless, this method does not enable to eliminate  the noise from the measurements.
See appendix \ref{mat_lissageannexe}.
{The smoothing parameter is automatically determined by minimizing and validated by figure \ref{reactionXY}.}

\subsection{Inverse dynamics   method and methods to determine $I_1$, $I_2$, $I_3$ and  $I_4$}
\label{mat_DIF}

The dynamics  equations applied to each of the segments $[A_jA_{j+1}]$, for $j\in \{1,...,q-1\}$,  give
\begin{subequations}
\label{RFDtot}
\begin{align}
\label{40}
&\vec R_j-\vec R_{j+1}-=-m_j\vec g+m_j\frac{d^2 \overrightarrow{OG_j}}{dt^2},\\
\label{60}
&-{\mathcal{M}}_j+I_j\ddot \phi_j=   C_j-C_{j+1},
\end{align}
\end{subequations}
where 
\begin{equation}
\label{60.a}
{\mathcal{M}}_j=-(x_{j+1}-x_j)\left(\alpha_jR_{y,j}+(1-\alpha_j)R_{y,j+1}\right)
+(y_{j+1}-y_j)\left(\alpha_jR_{x,j}+(1-\alpha_j)R_{x,j+1}\right).
\end{equation}
With boundary condition 
\begin{align}
\label{45tot}
&\vec R_1=\vec R,\quad 
\vec R_p=\vec 0,
\\
\label{70tot}
&C_1=C,\quad
C_q=0,
\end{align}
we obtain classically (see \cite{Hof1992}), for all $k\in\{1,...,q-1\}$, 
\begin{subequations}
\label{dynainvtot}
\begin{align}
\label{46tot}
& \vec R_{k}=\vec R-\sum_{j=1}^{k-1}  m_j\left(\frac{d^2{\overrightarrow{OG_j}}}{dt^2}-\vec g\right),\\
& C_k=C+\sum_{j=1}^{k-1}\left( {\mathcal{M}}_j-I_j\ddot \phi_j\right),
\intertext{and} 
\label{80}
&C=-\sum_{j=1}^{q-1}   {\mathcal{M}}_j+\sum_{j=1}^{q-1}I_j\ddot \phi_j.
\end{align}
\end{subequations}
The residual torque is defined by 
\begin{equation}
\label{90}
\widetilde C=C+\sum_{j=1}^{q-1}   {\mathcal{M}}_j-\sum_{j=1}^{q-1}I_j\ddot \phi_j,
\end{equation}
where angles 
$\phi_j$ are determined  from the smoothed displacements,   ${\mathcal{M}}_j$ 
are defined by  \eqref{60.a}   and joint forces   $R_{x,j}$ et $R_{y,j}$ are calculated  by using \eqref{46tot}.

We now explain how to determine $I_1$, $I_2$, $I_3$ and  $I_4$.

The residual torque is defined by \eqref{90} 
or
by the following equation:
\begin{subequations}
\label{coupleresi0}
\begin{align}
\label{coupleresiA0}
&\widetilde C^{(0)}(t)=C_{\text{exp}}-C_{\text{angl}},
\intertext{where $C_{\text{exp}}$ is torque measured
experimentally and}
\label{coupleresiB0}
&C_{\text{angl}}=-\sum_{j=1}^{q-1}   {\mathcal{M}}_j+
\sum_{j=1}^{q-1}I_j\ddot \phi_j,
\end{align}
\end{subequations}
is defined according to moments ${\mathcal{M}}_j$ and the double derivatives $\ddot \phi_j$.
$X^{(0)}$ corresponds to the values of function $X$.

The impulsion phase is equal to $[t_0,t_f]$.

By integration, between the beginning $t_0$  and $t_i$, we obtain,
and since the angular velocities
are null at the onset
of the push-off, we obtain 
\begin{subequations}
\label{coupleresi1}
\begin{align}
\label{coupleresiA1}
&\widetilde C^{(1)}(t_i)=C^{(1)}_{\text{exp}}(t_i)-C_{\text{angl}}^{(1)}(t_i),\\
\label{coupleresiAb1}
&C^{(1)}_{\text{exp}}(t_i)=\int_{t_0} ^{t_i} C_{\text{exp}}(s)ds,\\
\label{coupleresiB1}
&C_{\text{angl}}^{(1)}(t_i)=-\sum_{j=1}^{q-1}  \int_{t_0} ^{t_i}  {\mathcal{M}}_j(s)ds+
\sum_{j=1}^{q-1}I_j\dot \phi_j(t_i),
\end{align}
\end{subequations}
where $s$ is the variable of integration.
$X^{(1)}$ corresponds to the
first order integration of
the function $X$.
After a
second integration we obtain:
\begin{subequations}
\label{coupleresi2}
\begin{align}
\label{coupleresiA2}
&\widetilde C^{(2)}(t_i)=C_{\text{exp}}^{(2)}(t_i)-C_{\text{angl}}^{(2)}(t_i),\\
\label{coupleresiAb2}
&C^{(2)}_{\text{exp}}(t_i)=\int_{t_0} ^{t_i} \int_{t_0} ^uC_{\text{exp}}(s)dsdu,\\
\label{coupleresiB2}
&C_{\text{angl}}^{(2)}(t_i)=-\sum_{j=1}^{q-1}  \int_{t_0} ^{t_i} \int_{t_0} ^{u}  {\mathcal{M}}_j(s)dsdu+
\sum_{j=1}^{q-1}I_j\left(\phi_j(t_i)-\phi_j(t_0)\right),
\end{align}
\end{subequations}
where $u$ is the second variable of integration.
$X^{(2)}$ corresponds to the
second order integration of
the function $X$.
In order to compare the residual values $\widetilde C^{(0)}$, $\widetilde C^{(1)}$ and $\widetilde C^{(2)}$ obtained with  different methods,
it is necessary to normalize these values by considering  the dimensionless quantity defined by 
\begin{equation}
\label{100}
\varepsilon^{(j)}=\frac{\vnorm{C_{\text{exp}}^{(j)}-C^{(j)}_{\text{ang}}}}{\vnorm{C^{(j)}_{\text{exp}}}+\vnorm{C^{(j)}_{\text{ang}}}}\in [0,1],
\end{equation}
where $\vnorm{}$ is the $l^2$ norm, defined by \eqref{27}.


\subsubsection{Method A: optimization on all inertia $I_1$, $I_2$, $I_3$ and  $I_4$}\
\label{methA}

Considering that the residual is null, \eqref{coupleresi0},  \eqref{coupleresi1}, and  \eqref{coupleresi2} become 
\begin{subequations}
\label{90tot}
\begin{align}
\label{90A}
& \sum_{j=1}^{q-1}I_j\ddot \phi_j(t_i)=C(t_i)+\sum_{j=1}^{q-1}   {\mathcal{M}}_j(t_i),
\intertext{or}
\label{90B}
&\sum_{j=1}^{q-1}I_j\dot \phi_j(t_i)=\int_{t_0}^{t_i}  \left(C(s)+\sum_{j=1}^{q-1}   {\mathcal{M}}_j(s) \right)\, ds,
\intertext{or}
\label{90C}
& 
\sum_{j=1}^{q-1}I_j(\phi_j(t_i)-\phi_j(t_0))=\int_{t_0}^{t_i} \int_{t_0}^u \left(C(s)+\sum_{j=1}^{q-1}   {\mathcal{M}}_j(s) \right)\, dsdu.
\end{align}
\end{subequations}

As the method used in Section \ref{syntrraireac} to determine $\alpha_4$, for
Eq. \eqref{coupleresi2},  the double derivative of angles is not used for \eqref{90C}, but only values of these angles.

Each  equation  \eqref{90tot}    is equivalent to  determine $I_1$, $I_2$, $I_3$ and  $I_4$
such that 
\begin{equation}
\label{91}
\forall i, \quad \sum_{j=1}^{q-1}A_{i,j} I_j=B_i
\end{equation}
where $A_{i,j}$ and $B_i$ are known.
Theses equations are equivalent to the overdetermined linear system
\begin{equation}
\label{92}
A I=B,
\text{ where } 
I=
\begin{pmatrix}
I_1\\I_2\\I_3\\I_4
\end{pmatrix}
\end{equation}
which has no solution in the general case, but has a least square sens solution:
See appendix \ref{llsqannexe}.
In this case, 
the number $j\in \{0,1,2\}$ is called the degree of the method A; 
the number $\varepsilon^{(j)}$, defined by \eqref{100}
is denoted $\varepsilon^{(j)}_\text{A}$
and the coefficient of multiple determination for the overdetermined  system  
\eqref{92} is denoted ${R^2}^{(j)}_\text{A}$.

\subsubsection{Method B: optimization only on inertia $I_4$}\
\label{methB}

It can be   assumed that $I_1$, $I_2$, and $I_3$ are determined in \cite{winter}.
Then \eqref{coupleresi0}, \eqref{coupleresi1}, or \eqref{coupleresi2}  can be written under the following form: for all $i$, 
\begin{subequations}
\label{99Atot}
\begin{align}
\label{99Aa}
& I_{q-1}\ddot \phi_{q-1}(t_i)=-\sum_{j=1}^{q-2}I_{j}\ddot \phi_{j}(t_i)+C(t_i)+\sum_{j=1}^{q-1}   {\mathcal{M}}_j(t_i),\\
\label{99Ab}
&I_{q-1}\dot \phi_{q-1}(t_i)=- \sum_{j=1}^{q-2}I_j\dot \phi_j(t_i)+\int_{t_0}^{t_i}  \left(C(s)+\sum_{j=1}^{q-1}   {\mathcal{M}}_j(s) \right)\, ds,
\intertext{or} 
& \label{99Ac}
I_{q-1}(\phi(t_i)-\phi_{q-1}(t_0))=- \sum_{j=1}^{q-2}I_j(\phi_j(t_i)-\phi_j(t_0))+\int_{t_0}^{t_i} \int_{t_0}^u \left(C(s)+\sum_{j=1}^{q-1}   {\mathcal{M}}_j(s) \right)\, dsdu.
\end{align}
\end{subequations}

Each least square linear \eqref{99Atot} system can be written under the form
\eqref{92.A}-\eqref{92.B} of appendix \ref{llsqannexe}. Here, it is also equivalent to: find $I_4$
such that 
\begin{equation}
\label{99B}
\forall i,\quad 
y_i=I_4x_i.
\end{equation}
As previously,  
we consider  $\varepsilon^{(j)}_\text{B}$  and ${R^2}^{(j)}_\text{B}$.

\subsubsection{Method C: values of  $I_1$, $I_2$, $I_3$ and  $I_4$ defined by Winter}\
\label{methC}

The values of $I_1$, $I_2$, $I_3$ and  $I_4$ are estimated from \cite{winter}.
As previously, 
we consider  $\varepsilon^{(j)}_\text{C}$  and ${R^2}^{(j)}_\text{C}$.
This method is not an optimization method and ${R^2}^{(j)}_\text{C}$ is formally defined;
this number is not necessarily positive.

\vspace{1 cm}

To summarize,  we have three  methods defined by $\text{X} \in \{\text{A},\text{B},\text{C}\}$ 
and for each of them the order $j$ belongs to $j\in \{0,1,2\}$. 
The method $X$ with degree $j$ is called method "X$j$". For example "A2" is the method A with degree 2.
For each of these three methods and for each degree $j$
are defined $\varepsilon^{(j)}_\text{X}$ and ${R^2}^{(j)}_\text{X}$.
An accurate method corresponds to $\varepsilon$, close to 0 and $R^2$ close to 1.

\subsection{Statistics}
\label{dhdhs}

Main effects of the three methods and the three degrees based  on "residual error" were tested to significance with a general linear model one way ANOVA 
for repeated measures. When a significant F value was found, post-hoc Tukey tests were applied  to establish difference 
between methods (significant level $p< 0.05$) in Section \ref{generalisation}.

All analyses were proceeding through  the R software \cite{R}.

\section{Results}
\label{results}

The obtained results are relative to:
\begin{itemize}
\item[$\bullet$]
trunk anthropometry: 
\begin{itemize}
\item[-]
value of $\alpha_4$;
\item[-]
value of $I_4$;
\end{itemize}
\item[$\bullet$]
joint forces  $\vec R_k$ and torque $C_k$.
\item[$\bullet$]
and for each method 
by $\text{X} \in \{\text{A},\text{B},\text{C}\}$ 
and degree  $j\in \{0,1,2\}$ and called "X$j$":
\begin{itemize}
\item[-]
$\varepsilon^{(j)}_\text{X}$ defined by    \eqref{100};
\item[-]
${R^2}^{(j)}_\text{X}$, the coefficient of multiple determination, for the  system  
\eqref{92.A}-\eqref{92.B}.
\end{itemize}
\end{itemize}

\subsection{Validation of procedures and methods for one subject}
\label{validprocedure}

First, the validation will concern the synchronization, the smoothing, the three  methods with three degrees  
and the inverse dynamic method and will be illustrated for one subject.

\subsubsection{Synchronization}\
\label{validation_synchro}

\begin{figure}[ht]
\begin{center}
\epsfig{file=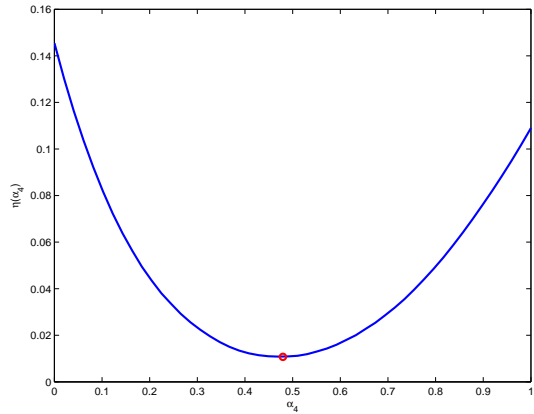,  width=7 cm}
\end{center}
\caption{\label{courbeeta} Curve $\eta(\alpha_4)$ {defined by \eqref{27c}}. The optimal $\alpha_4$ {which corresponds to 
the mimimum value of $\eta$ }is plotted by a red circle.}
\end{figure}

The figure \ref{courbeeta} presents the evolution of $\eta$ related to $\alpha_4\in [0,1]$. The
optimal value of $\alpha_4$ is plotted by a red circle; it  is given by 
\begin{equation}
\label{28}
\alpha_4=\alphaquatre,
\end{equation}
that can be compared to the value determined from  \cite{winter}
\begin{equation}
\label{eqalphawinter}
\alpha_4^\text{W}=\alphaquatrewinter.
\end{equation}
As assumed previously, this difference is related to the difference in the position 
 of HAT  during squat jumping and the position collected on cadavers.

\begin{figure}[ht]
\psfrag{Force (N)}[][l]{{\scriptsize{Vertical displacement of the body center of mass (m.)}}}
\begin{center}
\epsfig{file=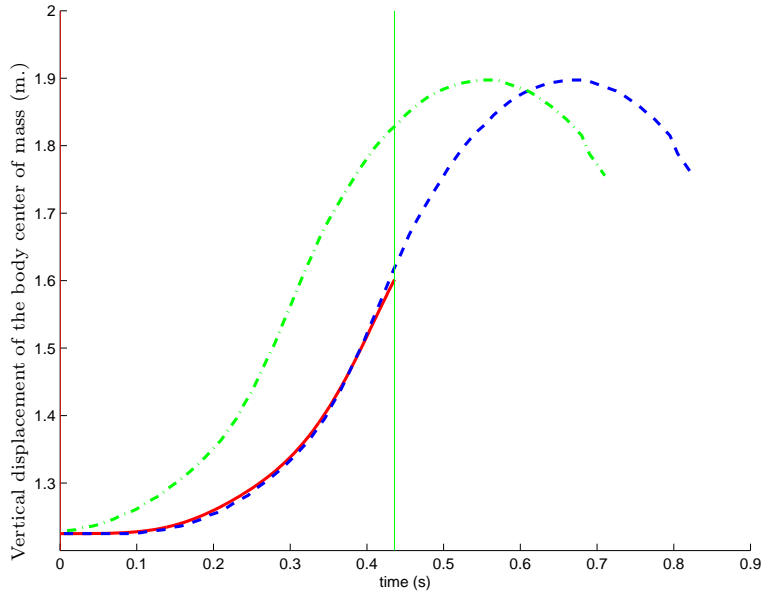,  width=10 cm}
\end{center}
\caption{\label{ordonneeCG} Three methods to determine the ordinate of the center of mass for one subject: the 
double integration of the experimental force (red continuous line),
with synchronized experimental data (blue dashed line), and with non synchronized experimental (dashdot green line).
The ordinate of the center of mass determined with synchronized experimental data depend on $\nu$ and $\alpha_4$.
A green vertical line shows the end of the push-off.}
\end{figure}

See also figure \ref{ordonneeCG}.
The 
ordinates  of $G$ determined 
by synchronized experimental data  are close to 
those  determined with the double integration  of the experimental force.
This is not the case for non synchronized experimental data.
{This finding validates the synchronization.}

\subsubsection{Smoothing}\
\label{validation_lissage}

\begin{figure}[ht]
\centering
\subfigure[\label{reactionX}: horizontal  ground reaction]
{\epsfig{file=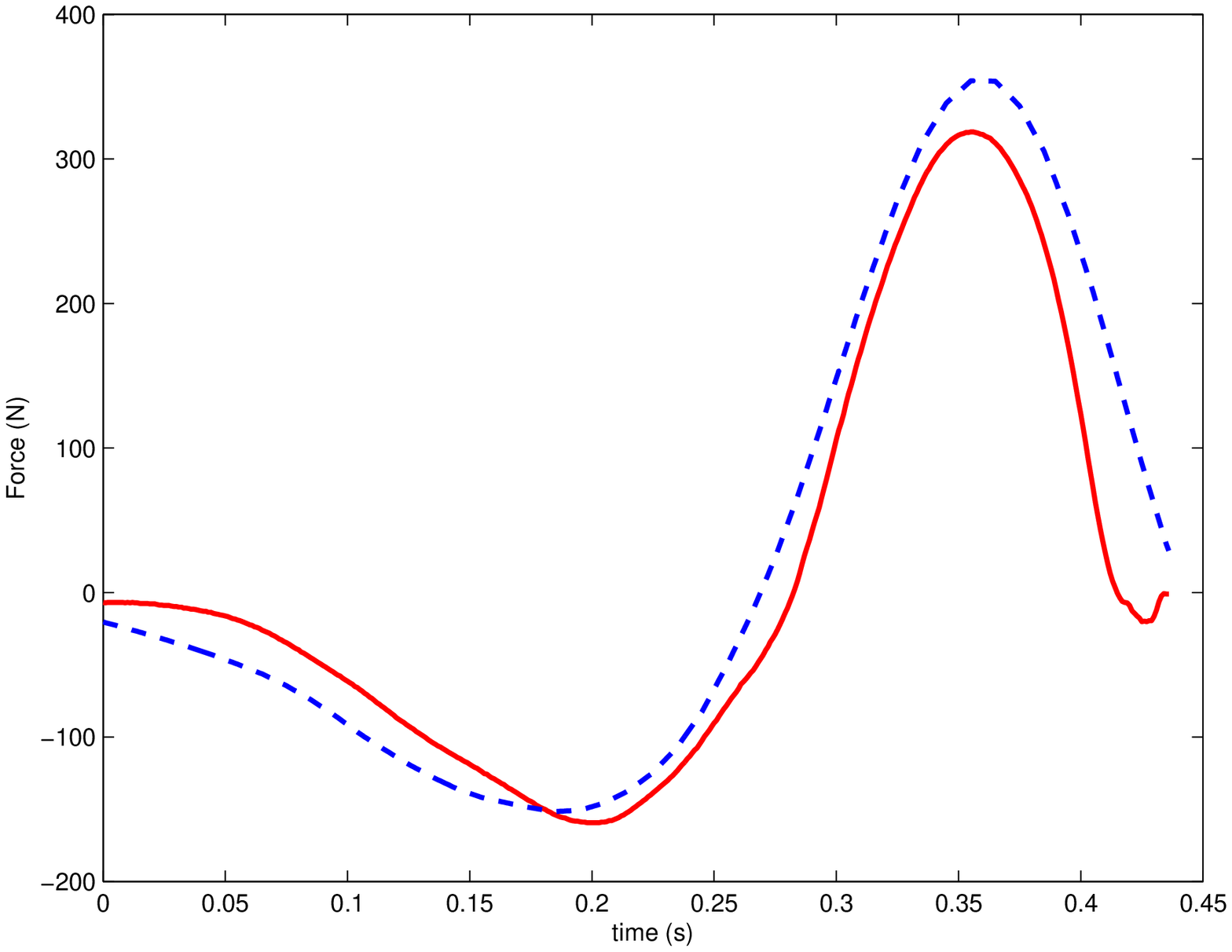, width=7 cm}}
\qquad
\subfigure[\label{reactionY}: vertical   of ground reaction]
{\epsfig{file=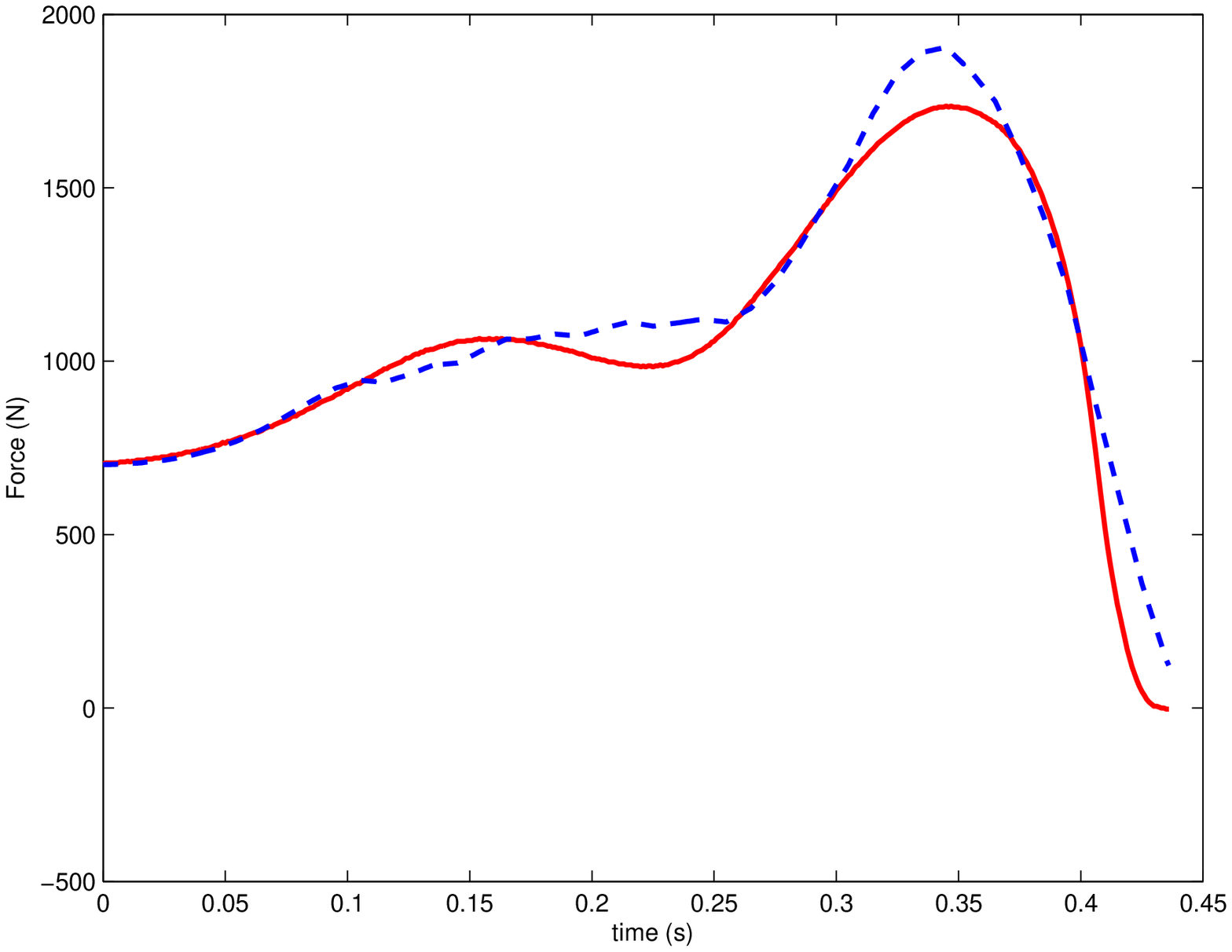, width=7 cm}}
\caption{\label{reactionXY} Ground reaction determined by two methods:
Experimental (red continuous line),
and by using a  smoothed acceleration of the center of mass  (blue dashed line).}
\end{figure}

Concerning the smoothing, the figure \ref{reactionXY} points out that curves of the experimental and calculated  reaction force
during the push-off phase, are close. From these results, the smoothing developped is considered as valuable.

\subsubsection{The three  methods and the three degrees}\
\label{validation_quatremetho}

\paragraph{Method A}\
\label{validation_methodA}

\begin{table}[t]
\caption{Values of obtained inertia $I_1$, $I_2$, $I_3$ and $I_4$ in  kgm$^2$  for  different cases.}
\begin{center}
\label{tab01}
\begin{tabular}{llll}
\hline
case & by values  (degree $j=0$) )& by integration  (degree  $j=1$)& by double integration  (degree $j=2$)\\
\hline
$I_1$ &\inertieunval &\inertieunint &\inertieunintd \\
$I_2$ &\inertiedeuxval &\inertiedeuxint &\inertiedeuxintd \\ 
$I_3$ &\inertietroisval & \inertietroisint& \inertietroisintd\\ 
$I_4$ &\inertiequatreval & \inertiequatreint&\inertiequatreintd \\ 
\hline
\end{tabular}
\end{center}
\end{table}

\begin{table}[t]
\caption{Values of $\alpha_i$,$m_i/m$,  $\widetilde r_i$, and $I_i$ (in  kgm$^2$) according to  \cite{winter}.}
\begin{center}
\label{tab00}
\begin{tabular}{lllll}
\hline
$i$ & $\alpha_i$ &  $m_i/m$ & $\widetilde r_i$ &  $I_i$ \\
\hline
 $1$ & \alphaunwinter & \mparunwinter & \rnormunwinter & \inertieunwinter\\ 
 $2$ & \alphadeuxwinter & \mpardeuxwinter & \rnormdeuxwinter & \inertiedeuxwinter\\ 
 $3$ & \alphatroiswinter & \mpartroiswinter & \rnormtroiswinter & \inertietroiswinter\\ 
 $4$ & \alphaquatrewinter & \mparquatrewinter & \rnormquatrewinter & \inertiequatrewinter\\ 
\hline
\end{tabular}
\end{center}
\end{table}

The calculated IP are presented in Table \ref{tab01}.
This method is not valuable while it provides values without any physical meaning. Indeed, we can see that:
\begin{enumerate}
\item
For $j=0$, obtained values are note positive.

In this case, we may consider a constrained optimization with inequalities $I_i\geq 0$.
However, 
since the unconstrained overdetermined linear problem possesses a unique solution, the  constrained optimization  problem will give solution 
where at least one inertia is equal to zero, which is not interesting  from a physical viewpoint.

\item
For $j=1$ or $j=2$, obtained values differ greatly from values given by \cite{winter} from table \ref{tab00}.
\end{enumerate}
Then, this method has to be removed.

\paragraph{Method B}\
\label{validation_methodB}

\begin{figure}[ht]
\centering
\subfigure[\label{nuage0}: degree $j=0$]
{\epsfig{file=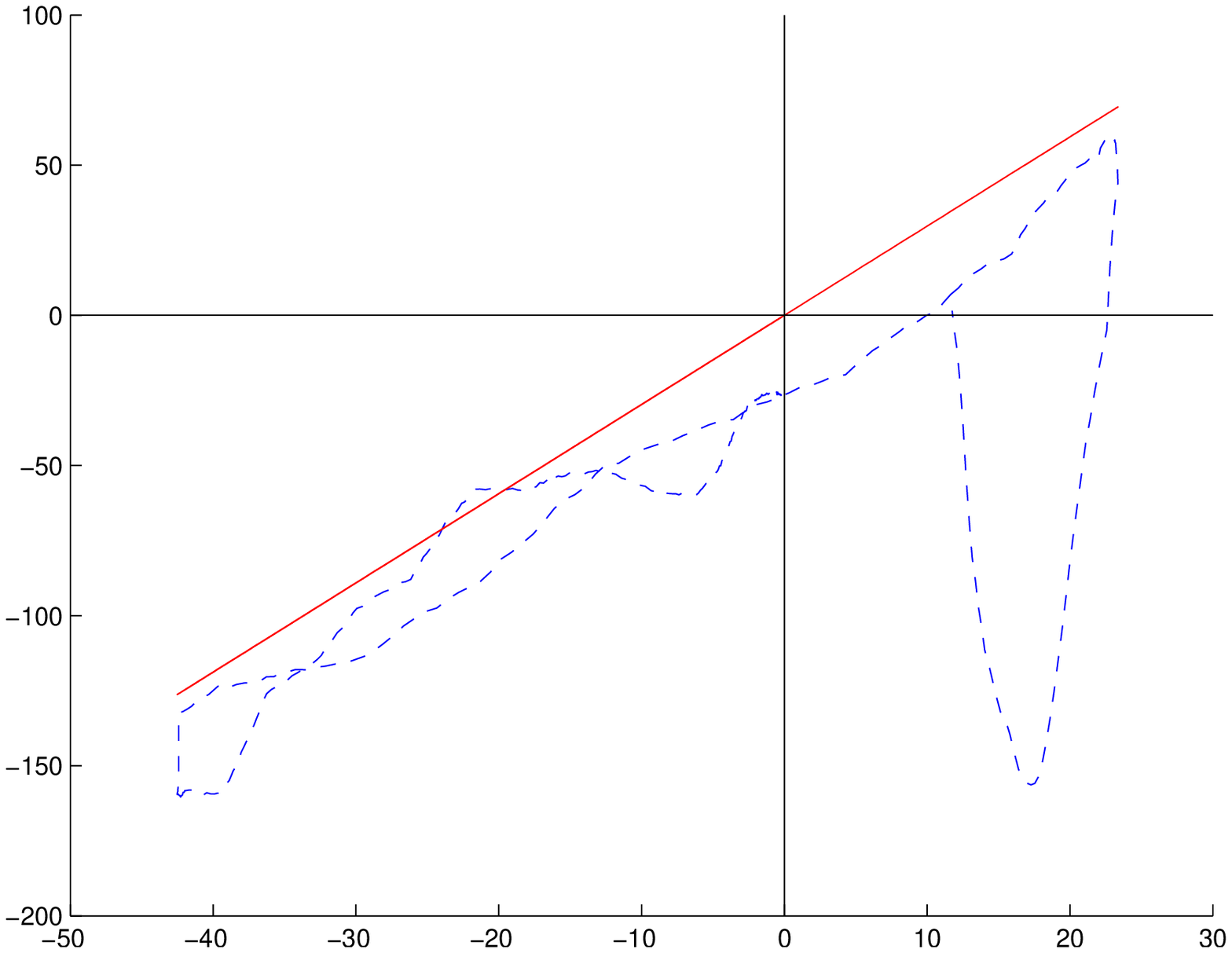, width=7 cm}}
\qquad 
\subfigure[\label{nuage1}:  degree $j=1$]
{\epsfig{file=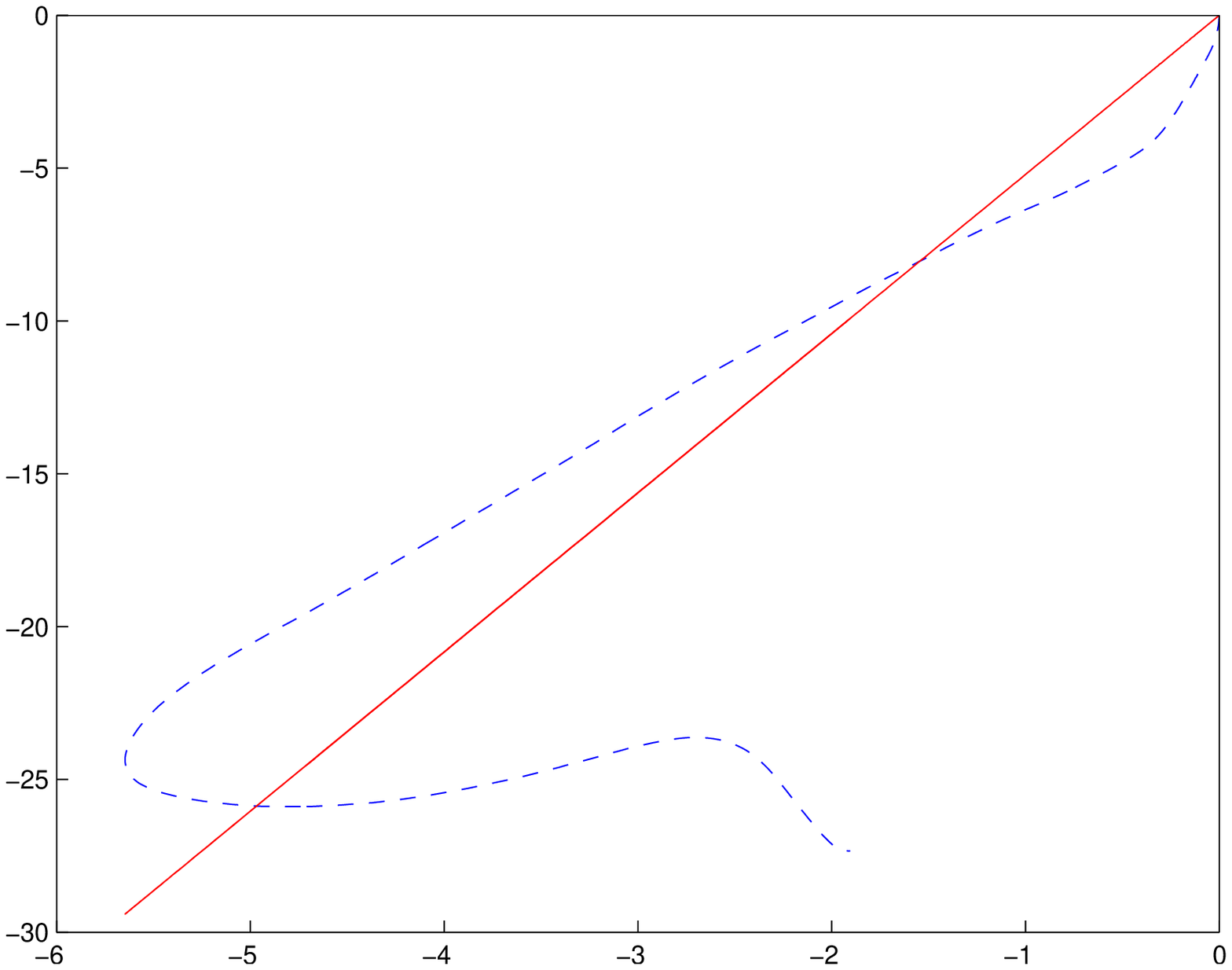, width=7 cm}}
\qquad 
\subfigure[\label{nuage2}:  degree $j=2$] 
{\epsfig{file=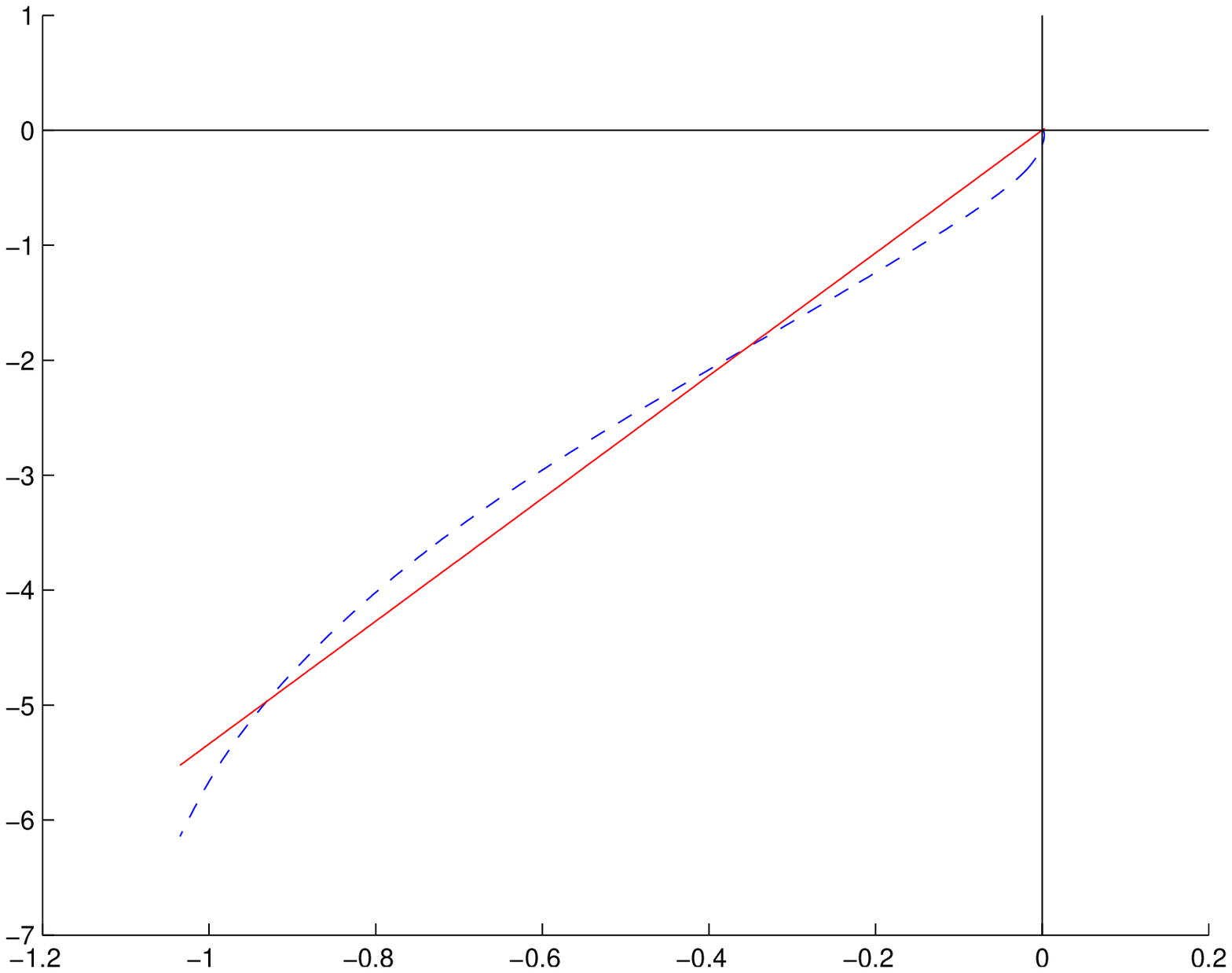, width=7 cm}}
\caption{\label{nuage} Points $(x_i,y_i )$  for different degrees;
point $(x_i,y_i)$ are plotted with  a blue  dashed line
and points $(x_i,I_4 x_i)$ are plotted with a red continuous line.}
\end{figure}
From the plotting  points $(x_i,y_i)$ defined by \eqref{99Atot} and \eqref{99B} for the three values of $j$
(in figure \ref{nuage}), 
it can be seen  that the best result  corresponds to  $j=2$. 
{%
For the Figure \ref{nuage2}, the slope is equal to 
$5.32608$. This value is greater to the one obtained from  Winter data ($\inertiequatrewinter$).}
The corresponding value of $\widetilde r_4$ is equal to 
\begin{equation}
\label{rquatrenorm}
\widetilde r_4=\rnormquatre,
\end{equation}
it can be compared to the value determined from \cite{winter}:
\begin{equation}
\label{eqrquatrenormwinter}
\widetilde r_4^\text{W}=\rnormquatrewinter.
\end{equation}

\paragraph{Method C}\
\label{validation_methodC_D}

For this method, 
inertia are chosen according to Winter (method C).

\paragraph{Comparison bettwen the three  methods and the three degrees and conclusion}\
\label{validation_methods}

\begin{table}[t]
\caption{Values of  $\varepsilon^{(j)}_\text{X}$ for different cases.}
\begin{center}
\label{tab11}
\begin{tabular}{llll}
\hline
degree $j$ & method A & method B & method C \\
\hline
0 & \residuvalA & \residuvalB & \residuvalC \\
1 & \residuintA & \residuintB & \residuintC \\
2 & \residuintdA & \residuintdB & \residuintdC \\
\hline
\end{tabular}
\end{center}
\end{table}

\begin{table}[t]
\caption{Values of ${R^2}^{(j)}_\text{X}$    for different cases.}
\begin{center}
\label{tab10}
\begin{tabular}{llll}
\hline
degree $j$ & method A & method B & method C \\
\hline
0 & \coeffdetvalA & \coeffdetvalB & \coeffdetvalC \\
1 & \coeffdetintA & \coeffdetintB & \coeffdetintC \\
2 & \coeffdetintdA & \coeffdetintdB & \coeffdetintdC \\
\hline
\end{tabular}
\end{center}
\end{table}


First of all, for $\varepsilon^{(j)}_\text{X}$ (table \ref{tab11}), 
the lower the value {close to 0}, the more accurate the method. 
Concerning
${R^2}^{(j)}_\text{X}$ (table \ref{tab10}), the more the value is close to 1, the more accurate the method. Therefore,
whatever the degree, $\varepsilon^{(j)}_\text{X}$ and ${R^2}^{(j)}_\text{X}$  are less and less accurate from methods A to C.
Considering the degree of integration, for the three methods (A to C), the results are more and more precise  
from degree 0 to degree 2.
Finally, for the particular studied subject, the tables  \ref{tab11} and \ref{tab10}
enable to established that the more accurate method was A2, followed by B2 and A1
 then by B1. We also see that C2 is more accurate than C0. 
We note that under the form 
\begin{subequations}
\label{resultsunsujet}
\begin{align}
\label{resultsunsujeta}
&\text{A2} <\text{B2} <\text{A1}<\text{B1},\\
\label{resultsunsujetb}
&\text{B2} <\text{B1} <\text{B0},\\
\label{resultsunsujetc}
&\text{C2} <\text{C0},
\end{align}
\end{subequations}
where  "$<$" means "more accurate than".


The method A 
can not be applied   regarding the   values given for inertia
$I_i$ (see section \ref{validation_methodA}).
Therefore, the best method, physically acceptable 
is the method B2.

This result will be confirmed by the statistics of section \ref{statscompmeth}.

\subsubsection{Inverse dynamic}\
\label{validation_inversemeth}

\begin{figure}[ht]
\centering
\subfigure[\label{RX}: horizontal joint reaction]
{\epsfig{file=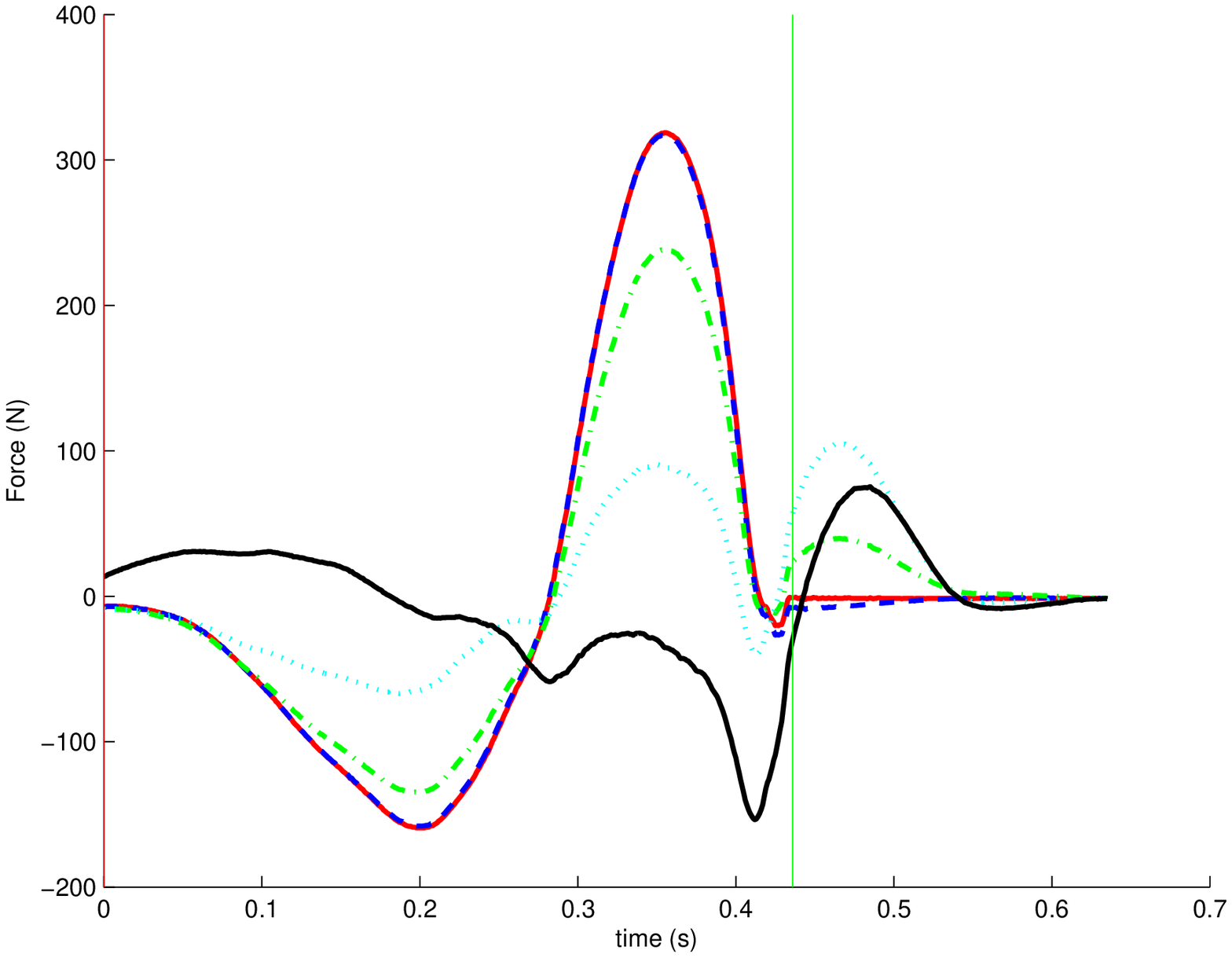, width=12 cm}}
\qquad 
\subfigure[\label{RY}: vertical joint reaction]
{\epsfig{file=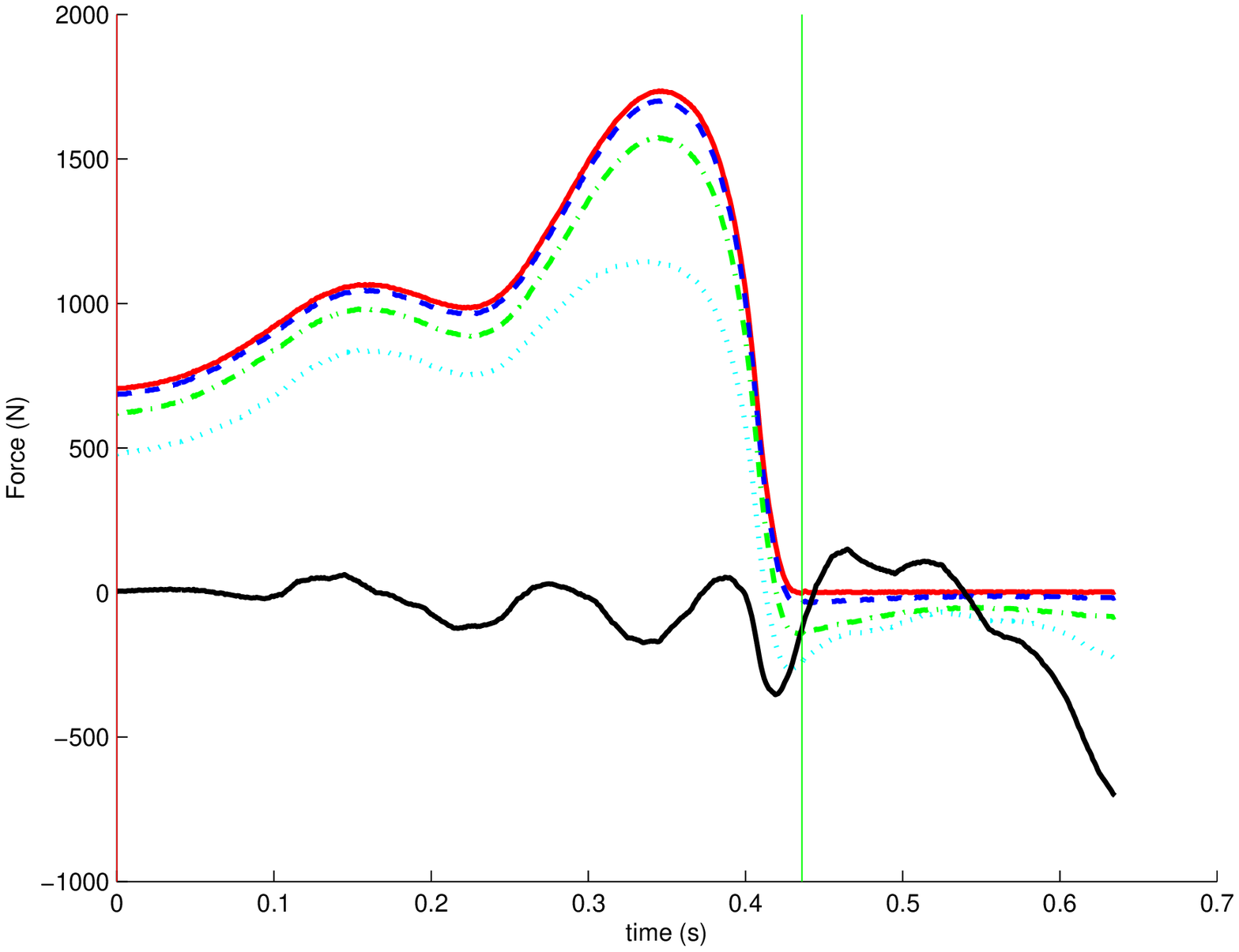, width=12 cm}}
\caption{\label{RXY} Differents joint reactions;
action from the support on the toe  are plotted by a red continuous line,
action from  the toe on the ankle by a blue dashed line,
action from  the ankle on knee by dashdot green line,
action knee  the on the hip by cyan dotted line, 
and residual action by a black continuous line.
The vertical green line corresponds to take off.
These values was determined from inertia given  by method B2.%
}
\end{figure}

\begin{figure}[ht]
\centering
\subfigure[\label{CX}: joint torque]
{\epsfig{file=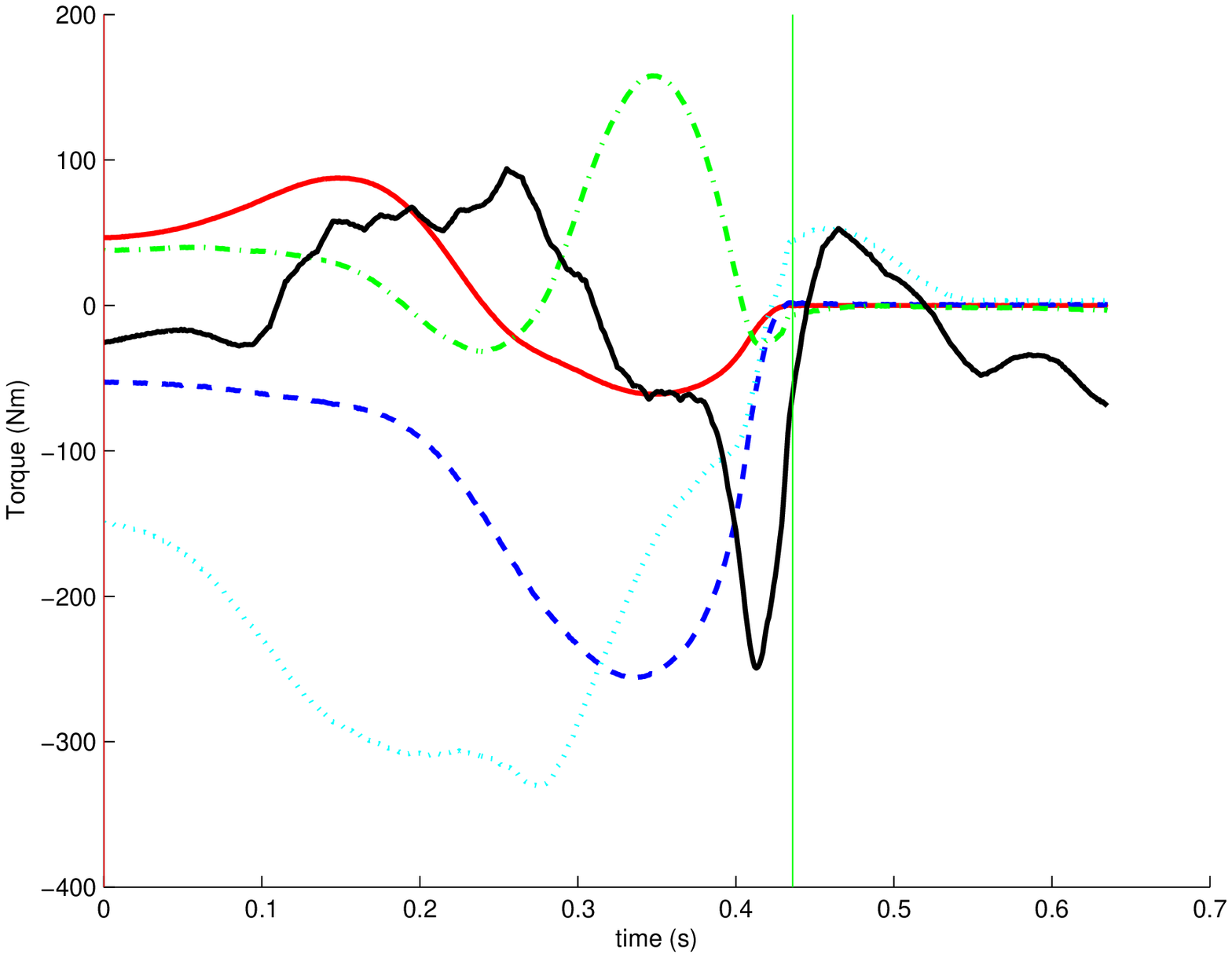, width=12 cm}}
\qquad 
\subfigure[\label{intdCX}: double integration of joint torque]
{\epsfig{file=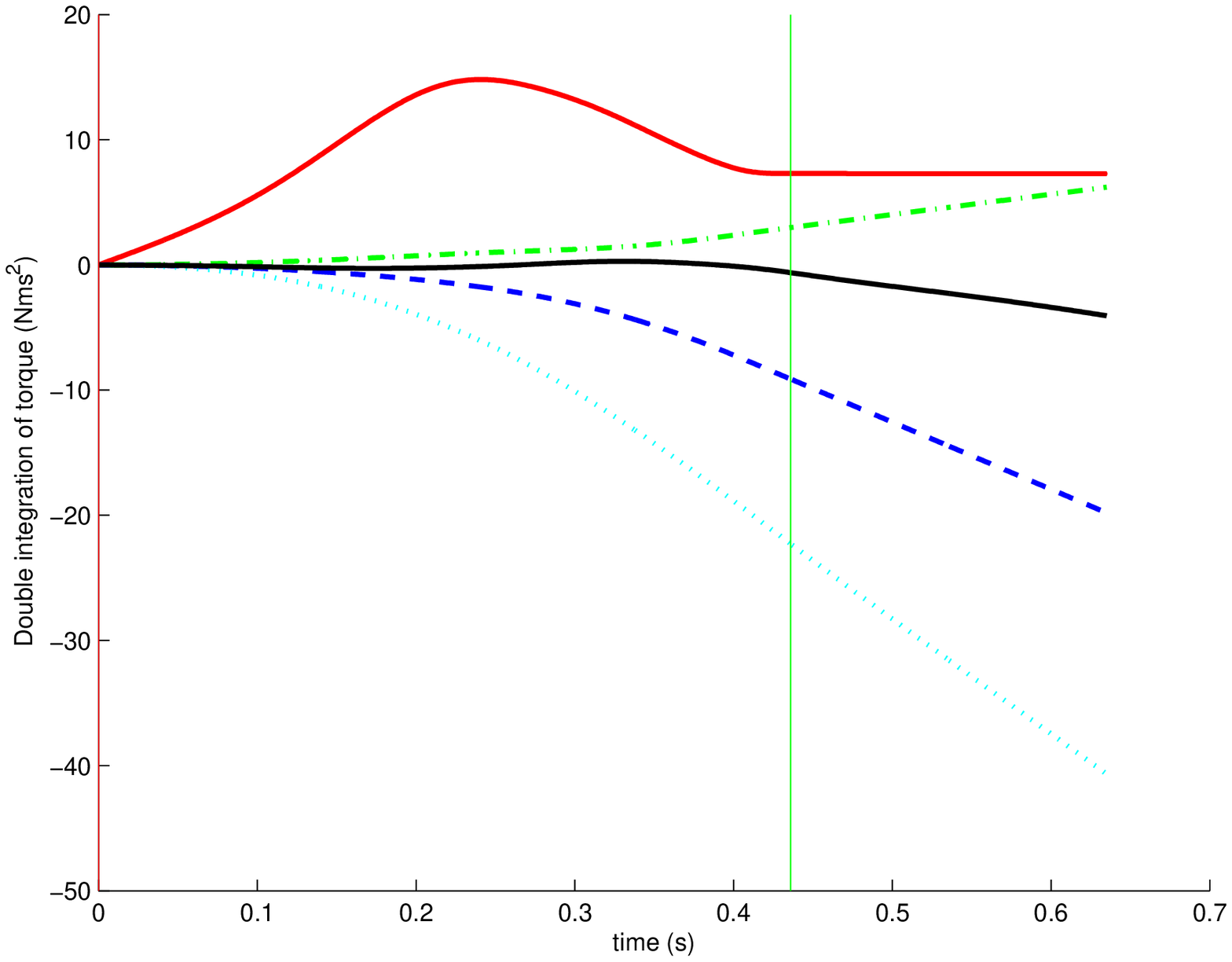, width=12 cm}}
\caption{\label{CXintdCX} Differents joint torques;
action from the support on the toe  are plotted by a red continuous line,
action from  the toe on the ankle by a blue dashed line,
action from  the ankle on knee by dashdot green line,
action knee  the on the hip by cyan dotted line, 
and residual action by a black continuous line.
The vertical green line corresponds to take off.
These values was determined from inertia given  by method B2.%
}
\end{figure}

The results of inverse dynamic methods are plotted on 
figures 
\ref{RXY}
 and  \ref{CXintdCX} for inertia given  by method B2.
For all figures, the beginning of the impulsion (corresponding to $t_0$) is plotted
by a vertical red line and 
the end  of the impulsion (corresponding to $t_f$) is plotted
by a vertical green line.

On figures \ref{CXintdCX},  the residual torque and their  double integrations 
are plotted. 
Torque value increased 
just before the end of impulsion. 
It can be noticed that  the norm $l^2$ of the double integral of this torque
were optimized but  
the maximum value of the torque were not optimized.
See figure   \ref{intdCX}: the smallest value of the double integrations corresponds to the residual torque.

\subsection{Generalization on the population (twelve subjets)}
\label{generalisation}




\setkeys{Gin}{width=0.55\textwidth}


12 subjects performed beetwen 
5 and 10  (mean: 7.25)
maximal squat jumps,
which gave   97 squat jumps.
The non positive or greater than 1 values of radius of gyration
were removed. Then, the number kept  for analysis was  87.

\subsubsection{Study of  $\alpha_4$ and $r_4$}
\label{statsalpha}

We now study $\alpha_4$
for trunk.

\begin{figure}[ht]
\begin{center}
\includegraphics{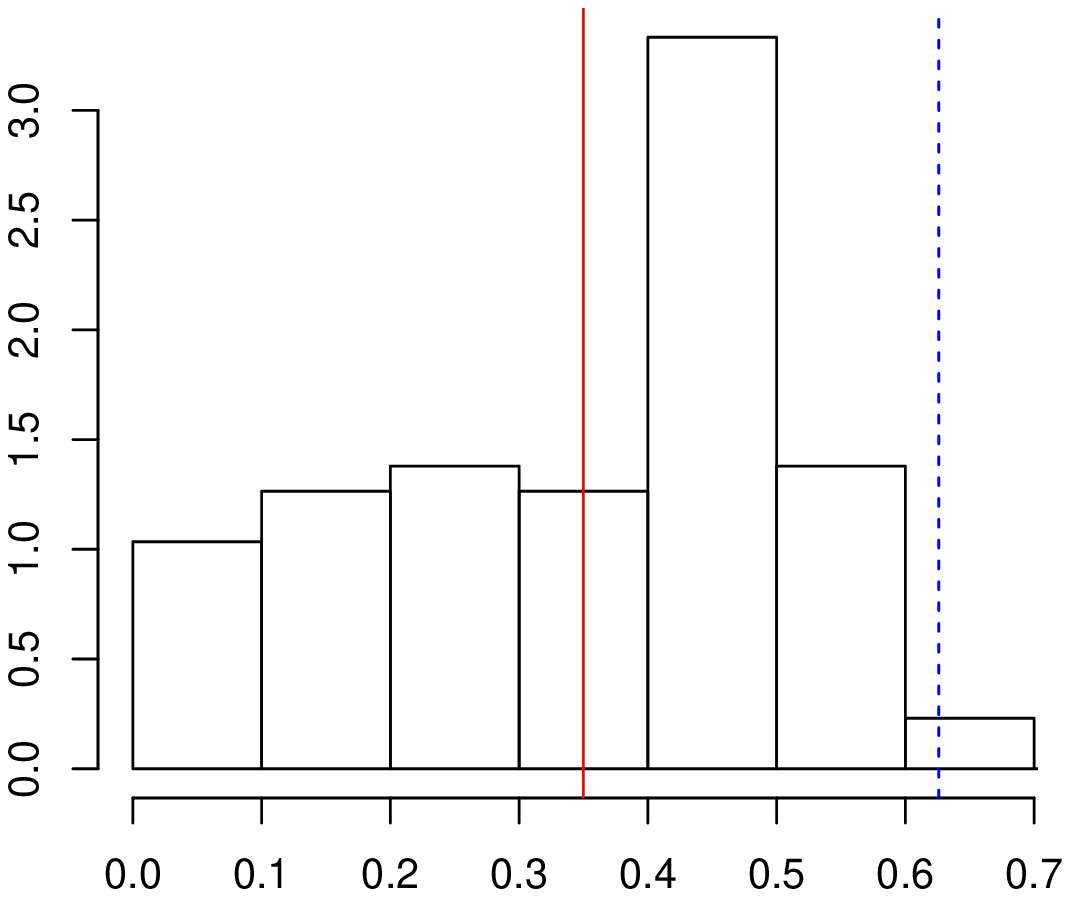}
\end{center}
\caption{\label{statsfig10} Histogram of $\alpha_4$.}
\end{figure}

See histogram in figure \ref{statsfig10}.
On all figures, we added the mean of data, plotted by a red continuous line, 
and  the value determined by Winter,  plotted by a blue dashed  line.
Here, according to Winter, $\alpha_4^\text{W}$ is defined by \eqref{eqalphawinter}.

\begin{table}[ht]
\caption{Basic statistics for $\alpha_4$ and $\widetilde r_4$}
\begin{center}
\label{bassttasalpha}
\begin{tabular}{llllllll}
\hline
variable & mean & sd  & $25\%$ quantile &$75\%$ quantile &$2.5\% (Q_1)$ quantile &$97.5\% (Q_3)$ quantile\\
\hline
$\alpha_4$ &   0.3499 &0.1858 & 
0.2124 &0.4718 &0.0011 &0.6063     \\
$\widetilde r_4$ &   0.5838 &0.2206 &
0.4348 &0.7557 &0.1753 &0.9669     \\ 
\hline
\end{tabular}
\end{center}
\end{table}

Basic statistics are given in Table \ref{bassttasalpha}.
We see that $ 50 \%$ of values belong to the interval
$[0.2124,0.4718]$
and that $ 95 \%$ of values belong to the interval
$[0.0011,0.6063]$.

\begin{figure}[ht]
\begin{center}
\includegraphics{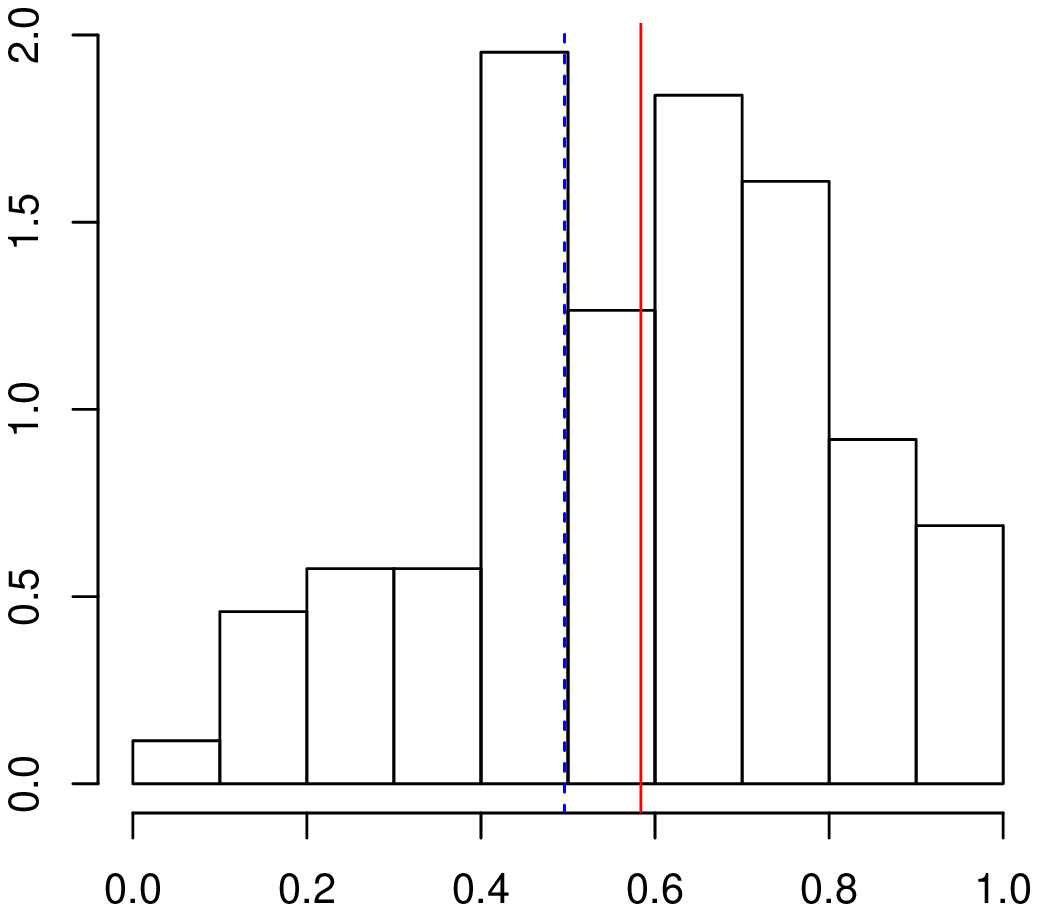}
\end{center}
\caption{\label{statsfig30} Histogram of $\widetilde r_4$.}
\end{figure}

We now study the normalized radius of gyration $\widetilde r_4$.
Here, for Winter, $\widetilde r_4^\text{W}$ is defined by \eqref{eqrquatrenormwinter}.
See figure \ref{statsfig30}.

Basic statistics are given in Table \ref{bassttasalpha}.
We see that $ 50 \%$ of values are included in the interval
$[0.4348,0.7557]$ and 
that $ 95 \%$ of values are included in the interval
$[0.1753,0.9669]$.

\subsubsection{Comparison beetwen Methods defined by X $\in$ $\{$A,B,C$\}$ and $j\in \{0,1,2\}$}
\label{statscompmeth}


The Shapiro-Wilk test shows that data $\varepsilon_\text{X}^{(j)}$ and  $1-{R^2}^{(j)}_\text{X}$
do not present a normal distribution; on the contrary, the logarithm of theses data follow a Gaussian distribution.

\begin{table}[ht]
\caption{Groups statistics of   $\log_{10}(\varepsilon_\text{X}^{(j)})$ for the three studied methods
with the three degrees: mean  $\pm $ standard deviation.}
\begin{center}
\label{statgroupeepsi}
\begin{tabular}{llll}
\hline
degree $j$ & method  A & method B  & method C    \\
\hline  
0 
& $-0.46 \pm 0.16$ 
& $-0.35 \pm 0.14$ 
& $-0.28 \pm 0.1$ 
\\
1
& $-1.06 \pm 0.29$ 
& $-0.59 \pm 0.3$ 
& $-0.37 \pm 0.25$ 
\\
2 
& $-1.83 \pm 0.44$ 
& $-1.01 \pm 0.38$ 
& $-0.49 \pm 0.33$ 
\\
\hline
\end{tabular}
\end{center}
\end{table}

\begin{table}[ht]
\caption{Groups statistics of  $\log_{10}\left(1-{R^2}^{(j)}_\text{X}\right)$  for the three studied methods
with the  three degrees: mean  $\pm $ standard deviation.}
\begin{center}
\label{statgrouper2}
\begin{tabular}{llll}
degree $j$ & method  A & method B  & method C    \\
\hline  
0 
& $-0.27 \pm 0.28$ 
& $-0.05 \pm 0.23$ 
& $0.11 \pm 0.24$ 
\\
1
& $-1.47 \pm 0.48$ 
& $-0.52 \pm 0.31$ 
& $0.03 \pm 0.54$ 
\\
2 
& $-2.99 \pm 0.73$ 
& $-1.37 \pm 0.56$ 
& $-0.2 \pm 0.82$ 
\\
\hline
\end{tabular}
\end{center}
\end{table}

Recall that for $p\in [0,1]$
\begin{itemize}
\item 
'***' means $p<0.001$;
\item
'**' means $p<0.01$;
\item
'*' means $p<0.05$;
\item
'.' means $p<.1$.
\end{itemize}

Statistics on $\log_{10}(\varepsilon_\text{X}^{(j)})$
and $\log_{10}\left(1-{R^2}^{(j)}_\text{X}\right)$ are given in tables \ref{statgroupeepsi} and \ref{statgrouper2}.
All tables with numerical details are given in Appendix \ref{annexetablestatsanovatukey}.

First of all, taken into consideration the $\log_{10}$
of 
$\varepsilon$ and $1-R^2$,
the lowest the value (close to 0), the more accurate the method or the degree of integration.

Both for 
$\log_{10}(\varepsilon)$ and $\log_{10}\left( 1-R^2 \right)$,
the general linear model one way ANOVA for repeated measures pointed out significant differences between the three methods (A, B and  C) and the three degrees (0, 1 and 2) 
($p\leq 3.933e-233$ (***)).

Both for 
$\log_{10}(\varepsilon)$ and $\log_{10}\left( 1-R^2 \right)$
the post- hoc Tuckey tests indicated that for the methods A and B, the values decreased when the degree increased 
(degree 2 $<$ degree 1 $<$ degree 0, $p<1.856e-08$   (***)).
Concerning the method C, the results of $\log_{10}(\varepsilon)$ and $\log_{10}\left( 1-R^2 \right)$
were lower for degree 2 than degree 1 
($p<0.01663$   (*)).
No significant difference was observed between the degree 1 and degree 0.

Concerning the comparison of the methods, for the degree 1 and 2, the values 
of $\log_{10}(\varepsilon)$
the following order was observed: A $<$ B $<$ C 
(A$<$B$<$C: $p\leq 3.869e-07$ (***).
Concerning the values 
of  $\log_{10}\left( 1-R^2 \right)$
were the lowest for method A, then B, then C  
(A$<$B$<$C, $p\leq 1.205e-13$ (***)).
With regard to the degree 0, the method A was significantly lower than method B 
($p\leq 0.02942$ (*)).
No difference was observed between the other methods.

Finally, when all methods and degrees were compared together, the lowest values of
$\log_{10}(\varepsilon)$ and $\log_{10}\left( 1-R^2 \right)$
were observed for method A2. 
The later was significantly lower than method A1 and method B2. 
By comparing A1 and B2, we obtain $p\geq 0.2719$.
The latter were significantly lower than the method B1.
Moreover, if we compare B1 to C2 and C2 to C0, we obtain
a maximum value of $p$ value equal to 
$0.0181$  (*).
We denote all theses results under the following form:
\begin{subequations}
\label{resultsstats}
\begin{align}
\label{resultsstatsa}
&\text{A2} <\text{B2}=\text{A1}<\text{B1}<\text{C2}<\text{C0},\\
\label{resultsstatsb}
&\text{B2} <\text{B1}<\text{B0},\\
\label{resultsstatsc}
&\text{A0}<\text{B0},
\end{align}
\end{subequations}
that confirms the results for one subject \eqref{resultsunsujet}.
As shown in the table \ref{tab01}, the results of method A gave unphysical values. Therefore the most accurate method was the method B2.

C0 and C2 are  the methods with torque value or double integration of torque
corresponding to Winter's data respectively.  A0 and B0 are optimization methods on
all IP or only trunk IP with Winter's data. 
These two  methods would correspond to the choice of methods used   by
\cite{Riemer2008,Riemer2009,Dumas2007}.

Moreover, the  geometrical mean of error method A2 is 
$6.5338$ smaller than 
the one of B2, what we denote under the form:
$
\text{A2}/\text{B2} \leq 6.5338
$.
See appendix \ref{amelirationerreur} for complete results.
We obtain the corroboration of \eqref{resultsstats}:
\begin{subequations}
\label{resultsstatsrap}
\begin{align}
\label{resultsstatsrapa}
&\frac{\text{A2}}{\text{B2}} \leq  6.5338,\quad
\frac{\text{B2}}{\text{B1}} \leq  2.669,\quad
\frac{\text{B1}}{\text{C2}} \leq  1.2634,\\
\label{resultsstatsrapb}
&\frac{\text{C2}}{\text{C0}} \leq  1.5924,\quad
\frac{\text{B1}}{\text{B0}} \leq  1.7228,\quad
\frac{\text{A0}}{\text{B0}} \leq  1.272.
\end{align}
\end{subequations}
The error values decrease compared to the previous methods observed in the litterature.


\section{Discussion}
\label{discussion}

The purpose of this study was to adjust AP and IP of the human segments during squat jumping in order to minimize error in joint torque. 
The results indicated that the method A2 minimized the most the residual torque (ie. $\varepsilon$ and $1- R^2$) following by the methods B2 and A1 being more accurate than the method B1. Nevertheless, the method A yields unrealistic $I_j$, therefore the most accurate method  retained was the method B2. Consequently, the optimization focused especially on the HAT inertial and anthropometric parameters.
It seems to be possible to optimize AP and IP of one segment when the others are known, but the simultaneous optimization of three AP and IP segments seems to be difficult. 
According to \cite{Riemer2009} IP
and AP optimized for three segments can not be considered as true, while if two of
three segments are known, the IP and AP of the last segment can be calculated.

The IP and AP found with the method B2 are close to the Winter ones but gave better residual joint torque. These differences could be obviously explained by the different position between the subjects performing squat jumps and cadavers. Especially the position of the arms was different, influencing the IP and AP of the HAT segment. 
\cite{Riemer2008,Riemer2009} used optimization techniques to solve inverse dynamic by determining numerical angles which minimize the difference between the ground reaction force measured and the ground reaction force calculated. They found segmental angles which minimize an objective function under equality and inequality constraints, by taking into account the difference ground reaction force measured and the ground reaction force calculated. From these angles, joint torques are determined.  Our approach is different: techniques of \cite{Riemer2008,Riemer2009} to determine angle and torque were not applied. Only experimental displacements smoothed (see section  \ref{mat_lissage}) were considered. Then, with a direct inverse method, joint forces and torques were deduced.  Finally, an optimization is made on the residual torque and force to determine values of AP and IP. This optimization is very simple and fast, since it  is based on the least square linear method (see  Eq. \eqref{92.A} and \eqref{92.B}). 
It can be noticed that even if the residual torque or force are minimized, the error at each joint may be increased 
\cite{Kuo1998,Riemer2008,Riemer2009}.
However, in our study and for the method B, we only optimized AP and IP of the trunk, consequently the joint torque and forces at the hip, knee and ankle joints remained unchanged.
The major difference with classical way is the following point: it is possible to minimize residual torque, or its integral or its double integral. 
The best method corresponds to the minimization of the double integral, which do not use the double derivative of angle.
Equations \eqref{resultsstats} and \eqref{resultsstatsrap}
show that our methods
seem  to be better than classical one  according statistics or error results.

\section{Conclusion}
\label{conclusion}

The optimization of inertial and anthropometric parameters seems to be necessary when researchers use inverse dynamic methods. Indeed, cadaver data lead to errors in the calculi of the joint torques. These ones could be reduced by optimization methods.
The present method of optimization,    
{based on the double integration of residual},
 has been applied on the HAT segment but could also be applied on more segments. Therefore, further studies focusing on cutting the HAT segment into more segments (the pelvis and the rachis for example) will use this method to calculate the anthropometric and inertial parameters of the new segments composing the HAT.



\appendix

\section{A few theoretical reminders}
\label{annexe}

\subsection{Synchronisation of displacements and forces and determination of $\alpha_4$}
\label{syntrraireacannexe}

\newcommand{\diW}{\mathcal{R}}

We have 
\begin{equation}
\label{01}
\forall t\in[t_0,t_f],\quad
\vec R=-m\vec g+m\frac{d^2 \overrightarrow{OG}}{dt^2}.
\end{equation}
With this equation, experimental values $\vec R$
and the acceleration of $G$, can be compared.
Nevertheless, 
calcul of the acceleration of $G$ from a double derivation of the numerical experimental data 
leads to inaccuracy;
moreover, data are not synchronized.
Let introduce the residual force defined by 
\begin{equation}
\label{20tot}
\vec{\widetilde  R}=\vec R  +m \vec g-m\frac{d^2 \overrightarrow{OG}}{dt^2}.
\end{equation}
from a theoretical point of view, $\vec{\widetilde  R}$ 
should be equal to zero.
Experimentally, this residual force is not equal to zero.
It depends  to $\alpha_1$, $\alpha_2$, $\alpha_3$, which are known  (from table \ref{tab00})
and  $\alpha_4$, unknow.
It also depends  of the time phase  beetwen the forces and the displacements.
Then 
it depends on $\nu$, defined by \eqref{sjeqdecalage}.
Moreover, 
to avoid the determination of acceleration of $G$, the double integration of 
the residual force is applied:
\begin{equation}
\label{25}
\overrightarrow{\diW}(t)=\int_{t_0}^t \int_{t_0}^u (\vec R(s)+m\vec g)\,dsdu-m\left(\overrightarrow{OG}(t)-\overrightarrow{OG}(t_0)\right)
\end{equation}
The measures on  $x$-axis are smaller than of the $y$-axis. Then the ordinate
of the residual force will be considered 
\begin{equation}
\label{26}
\diW_y(t)=\int_{t_0}^t \int_{t_0}^u (R_{y}(s)-mg)\,dsdu-my_G(t)+my_G(t_0),
\end{equation}
where $y_G(t)$ is the ordinate of the center of mass of the body.
$\diW_y(t_i)$ are defined by 
\begin{equation}
\label{26bis}
\forall i,\quad
\diW_y(t_i)=\int_{t_0}^{t_i} \int_{t_0}^u (R_{y}(s)-mg)\,dsdu-my_G(t_i)+my_G(t_0).
\end{equation}
The double integral can be numerically calculated  from experimental 
$R_y^i$ and from  integer $\nu$ defined by  \eqref{sjeqdecalage} is the delay between force and displacement and $G(t_i)$ can be determined according to experimental data 
$y_j^i$ (ordinates of anatomical landmarks)
and 
${{\alpha_1}}$, 
${{\alpha_2}}$, 
${{\alpha_3}}$, which are known  (from table \ref{tab00})
and  $\alpha_4$, unknow.

If $y={\left(y_i\right)}_{1\leq i \leq P}$ is an element of $\Er^P$, we note by $l^2$ norm:
\begin{equation}
\label{27}
\vnorm{y}=\sqrt{\sum_{i=1}^P y_i^2}.
\end{equation}
Let $i_f$ be a last integer $i$ such that $t_i\leq t_f$. 
Then, the number 
\begin{equation}
\label{27b}
\vnorm{\diW_y}=\vnorm{{\left(\diW_y(t_i)\right)}_{0\leq i \leq i_f}},  
\end{equation}
depends on  $\alpha_4$ and $\nu$.
For each value of $\alpha_4$,  the value of $\nu$, which minimizes $\vnorm{\diW_y}$, noted
$\eta(\alpha_4)$, is determined:
\begin{equation}
\label{27c}
\forall \alpha_4\in [0,1],\quad 
\eta(\alpha_4)=\min_{\nu} \vnorm{\diW_y}
\end{equation} 
Secondly  the value of $\alpha_4$ which minimizes $\eta$ is calculated:
\begin{equation}
\label{A9bis}
\eta(\alpha_4)=\min_{\alpha\in [0,1]} \eta(\alpha).
\end{equation}

Since $\nu$ and $\alpha_4$ are determined, 
$W_y(t)$ is "small" and can be rewritten 
under the form 
\begin{equation}
\label{26ter}
\int_{t_0}^t \int_{t_0}^u \left(\frac{1}{m}R_{y}(s)-g\right)\,dsdu-y_G(t)+y_G(t_0)\approx 0
\end{equation}
which leads
to determine then the ordinate of $G$ 
from three methods:
\begin{itemize}
\item[$\bullet$]
with 
the 
double integration of the experimental force;
\item[$\bullet$]
with 
synchronized experimental data;
\item[$\bullet$]
with 
non synchronized experimental data.
\end{itemize}

If experimental forces and displacements are synchronized, 
the value of  $\nu$ is known and $\vnorm{W_y}$ depends only on $\alpha_4$.
In this case, 
\begin{equation*}
\diW_y(t_i)=A_i\alpha_4+B_i,
\end{equation*}
where $A_i$ and $B_i$ are known. The determination of $\alpha_4$
which minimizes $\vnorm{W_y}$ is obtained by solving 
\begin{equation*}
\forall i,\quad
A_i\alpha_4+B_i=0
\end{equation*}
in the least square sens: see appendix \ref{llsqannexe}.

\subsection{Linear least squares system}
\label{llsqannexe}

For integers $P,Q$ such that $P\geq Q$
for all matrix $A\in {\mathcal{M}}_{P,Q}(\Er)$, for all $B\in \Er^P$, 
overdetermined linear system is considered 
\begin{equation}
\label{92.A}
AX=B,
\end{equation}
in the following sens: find $x\in \Er^Q$ such that 
\begin{equation}
\label{92.B}
\vnorm{Ax-B}=\min_{X\in \Er^Q}\vnorm{AX-B}.
\end{equation}
There is a unique solution if  the rank of matrix $A$ is equal to $Q$ \cite{MR1231568}.
A system like \eqref{92.B}
is called  a  linear least squares system.
On a theoretical point of view, the  unique solution of 
\eqref{92.A} or \eqref{92.B} is given by 
\begin{equation}
\label{92.C}
x={\left({}^tAA\right)}^{-1}({}^tA)B,
\end{equation}
where ${}^tA$ is the the transpose of the matrix $A$. The matrix ${\left({}^tAA\right)}^{-1}({}^tA)$ is some times called
pseudoinverse of $A$.

\subsection{Smoothing of experimental data $x_j^i$,  $y_j^i$, $x_{G_j}^i$, $y_{G_j}^i$, 
$x_{G}^i$ and $y_{G}^i$}
\label{mat_lissageannexe}

Since the values of $x_j^i$,  $y_j^i$, $x_{G_j}^i$, $y_{G_j}^i$, 
$x_{G}^i$ and $y_{G}^i$, are experimental 
data, they are not necessarily smooth and can not be derivatived ones or twice.

Then the following smoothing is used: it 
returns the cubic smoothing spline for the
given data $\left({\mathcal{X}}_i=i/f, {\mathcal{Y}}_i\right)$ ($f$ is the acquisition frequency) and depending on the smoothing parameter $s_p\in [0,1]$.
This smoothing spline  $F$  minimizes
\begin{equation*}
s_p\sum_i \left({\mathcal{Y}}_i-F({\mathcal{X}}_i)\right)^2+(1-s_p)\int (F'')^2.
\end{equation*}
  For  $s_p=0$, the smoothing spline is the least-squares straight line fit to
  the data, while, at the other extreme, i.e., for  $s_p=1$, it is the "natural" 
  or variational cubic spline interpolant.  The transition region between 
  these two extremes is usually only a  rather small range of values for $s_p$
  and its location strongly depends on the data sites.
Smoothing values of data $\left({\mathcal{Y}}_i\right)$ are then noted 
$F_{s_p}({\mathcal{X}}_i)$. The smoothing derivatives of data can be obtained. 
Moreover, values can be determined for each value of time. Then, 
it can be assumed
that \eqref{sjeq01} and  \eqref{sjeqdecalage} hold for the 
smallest value of the frequency, now noted $f'_e$.
Since time phase $\nu$ is now determined, 
\eqref{sjeq01} and \eqref{sjeqdecalage}  can be replaced by 
\begin{subequations}
\label{dataexp}
\begin{align}
\forall j\in \{1,...,q\}, \quad \forall  i\in \{0,...,n'\},\quad
\label{dataexpa}
& x_j^i=x_j(i/f'_e),\\
\label{dataexpb}
& y_j^i=y_j(i/f'_e),\\
\label{dataexpc}
&\vec R^i=\vec R(i/f'_e),\\
\label{dataexpd}
&C^i=C(i/f'_e).
\end{align}
\end{subequations}
where $f'_e=1000$ Hz. is the common acquisition frequency (for force and displacement).
The values of displacement, velocity and acceleration of experimental data, 
$x_j^i$,  $y_j^i$, $x_{G_j}^i$, $y_{G_j}^i$, 
$x_{G}^i$ and $y_{G}^i$
are replaced by 
their smoothing values.

An important choice is the value of  each  smoothing parameter ${s_p}_j\in [0,1]$ (for each point $A_j$) . The values
of ${s_p}_j\in [0,1]$ 
which minimize the residual force defined by 
$\vnorm{{\widetilde  R_x}}$ and $\vnorm{{\widetilde  R_y}}$
defined from \eqref{20tot} were calculated.

The derivatives of angles defined by \eqref{eq01tot}  can be determined.

\section{Complete numerical statistical resultats}
\label{annexetablestats}


\newpage

\subsection{ANOVA and post-hoc Tukey tests}
\label{annexetablestatsanovatukey}

\begin{table}[ht]
\caption{ANOVA for repeated measures}
\begin{center}
\label{tabanovaglob}
\begin{tabular}{lll}
\hline
variable & $F$ value  &  $p(>F)$  \\
\hline
$\log_{10}(\varepsilon)$ &308.2546 &  3.933e-233  (***) \\
$\log_{10}(1-R^2)$       &390.9661 &  5.969e-264  (***) \\  
\hline
\end{tabular}
\end{center}
\end{table}

\begin{table}[ht]
\caption{Post-hoc Tukey tests on $\log_{10}(\varepsilon)$}
\begin{center}
\label{tukey012res2}
\begin{tabular}{lllllllll}
\hline
& 
method  A &&
method B  & &
method  C &
\\
&
Estimate &  $Pr(<t)$ &
Estimate &  $Pr(<t)$ &
Estimate &  $Pr(<t)$ 
\\ 
\hline
1 $\leq $ 0 &
$-0.6086$  & $0$ (***)&
$-0.2362$  & $1.856e-08$ (***)&
$-0.0898$  & $0.07801$ (.)
\\
2 $\leq $ 1 &
$-0.7647$  & $0$ (***)&
$-0.4263$  & $0$ (***)&
$-0.1123$  & $0.01663$ (*)
\\
\hline
\end{tabular}
\end{center}
\end{table}

\begin{table}[ht]
\caption{Post-hoc Tukey tests on $\log_{10}\left(1-{R^2}\right)$.}
\begin{center}
\label{tukey0121mcd}
\begin{tabular}{lllllllll}
\hline
& 
method  A &&
method B  & &
method  C &
\\
&
Estimate &  $Pr(<t)$ &
Estimate &  $Pr(<t)$ &
Estimate &  $Pr(<t)$ 
\\ 
\hline
1 $\leq $ 0 &
$-1.1917$  & $0$ (***)&
$-0.4731$  & $3.323e-10$ (***)&
$-0.0781$  & $0.6292$ ( )
\\
2 $\leq $ 1 &
$-1.5217$  & $0$ (***)&
$-0.8501$  & $0$ (***)&
$-0.2238$  & $0.006464$ (**)
\\
\hline
\end{tabular}
\end{center}
\end{table}

\begin{table}[ht]
\caption{Post-hoc Tukey tests on  $\log_{10}(\varepsilon)$.}
\begin{center}
\label{tukeyABCDres2}
\begin{tabular}{lllllll}
\hline 
degree $j$ 
& A $\leq$ B & 
& B $\leq$ C  
\\
&
Estimate &  $Pr(<t)$ &
Estimate &  $Pr(<t)$ 
\\
\hline
0 &
$-0.1045$  & $0.02942$ (*)&
$-0.0674$  & $0.2618$ ( )
\\
1 &
$-0.4768$  & $0$ (***)&
$-0.2138$  & $3.869e-07$ (***)
\\
2 &
$-0.8152$  & $0$ (***)&
$-0.5279$  & $0$ (***)
\\
\hline
\end{tabular}
\end{center}
\end{table}

\begin{table}[ht]
\caption{Post-hoc Tukey tests on  $\log_{10}\left(1-{R^2}\right)$.}
\begin{center}
\label{tukeyABCD1mcd}
\begin{tabular}{lllllll}
\hline 
degree $j$ 
& A $\leq$ B & 
& B $\leq$ C  
\\
&
Estimate &  $Pr(<t)$ &
Estimate &  $Pr(<t)$ 
\\
\hline
0 &
$-0.2261$  & $0.005771$ (**)&
$-0.1545$  & $0.09916$ (.)
\\
1 &
$-0.9448$  & $0$ (***)&
$-0.5494$  & $1.205e-13$ (***)
\\
2 &
$-1.6164$  & $0$ (***)&
$-1.1757$  & $0$ (***)
\\
\hline
\end{tabular}
\end{center}
\end{table}

\begin{table}[ht]
\caption{Post-hoc Tukey tests on  $\log_{10}(\varepsilon)$.}
\begin{center}
\label{tukeycroise01}
\begin{tabular}{lll}
\hline 
test &
 Estimate &  $Pr(<t)$ 
\\
\hline 
A1$\leq$ B2 &
$-0.0505$  & $0.2956$ ( )\\
B1$\leq$ C2 &
$-0.1015$  & $0.01805$ (*)\\
C2$\leq$ C0 &
$-0.2021$  & $1.001e-06$ (***)
\\
\hline
\end{tabular}
\end{center}
\end{table}

\begin{table}[ht]
\caption{Post-hoc Tukey tests on  $\log_{10}\left(1-{R^2}\right)$.}
\begin{center}
\label{tukeycroise02}
\begin{tabular}{lll}
\hline 
test &
 Estimate &  $Pr(<t)$ 
\\
\hline 
A1$\leq$ B2 &
$-0.0947$  & $0.2719$ ( )\\
B1$\leq$ C2 &
$-0.3256$  & $1.161e-05$ (***)\\
C2$\leq$ C0 &
$-0.302$  & $5.039e-05$ (***)
\\
\hline
\end{tabular}
\end{center}
\end{table}


\subsection{Decrease of error}
\label{amelirationerreur}

We try to compare
$\log_{10}(\varepsilon^{(j)}_\text{X})$ and 
$\log_{10}(\varepsilon^{(j')}_\text{X'})$.
If $Q=87$ is the number of measures
and if we write, for each measure $i$,
\begin{equation*}
{\eta}_i=\log_{10}\left(  {\varepsilon^{(j)}_\text{X}}_i \right)
\end{equation*}
we obtain the geometric mean according to the arithmetric mean 
\begin{equation*}
{\left(\prod_i   {\varepsilon^{(j)}_\text{X}}_i  \right)}^{1/Q}=10^{\left(\frac{1}{Q}\sum_i {\eta}_i\right)}=10^{\overline{\log_{10}(\varepsilon^{(j)}_\text{X})}}
\end{equation*}
Thus,  
\begin{equation}
\label{eqdud40}
\frac{%
{\left(\prod_i{\varepsilon^{(j')}_\text{X'}}_i   \right)}^{1/Q}
}%
{%
{\left(\prod_i{\varepsilon^{(j)}_\text{X}}_i   \right)}^{1/Q}
}=
10^{\left(-\overline{\log_{10}(\varepsilon^{(j)}_\text{X})}+\overline{\log_{10}(\varepsilon^{(j')}_\text{X'})}\right)}
.
\end{equation}
It implies that, in geometrical mean, if the difference of the logarithms is equal to $\Delta L$, 
the error is divided by  $10^{\Delta L}$.

By using table \ref{statgroupeepsi}, we obtain tables \ref{amelierationpardegre} and \ref{amelierationparmethode}

\begin{table}[ht]
\caption{Division of error according to method.}
\begin{center}
\label{amelierationpardegre}
\begin{tabular}{lll}
\hline 
degree $j$ &  
A/B &
B/C
\\
\hline 
0 
& 1.272
& 1.1678
\\
1
& 2.9979
& 1.6362
\\
2
& 6.5338
& 3.3719
\\
\hline
\end{tabular}
\end{center}
\end{table}

\begin{table}[ht]
\caption{Division of error according to degree.}
\begin{center}
\label{amelierationparmethode}
\begin{tabular}{llll}
\hline 
degree $j$ & method  A & method B  & method C
\\
\hline 
$1/0$
& 4.0605
& 1.7228
& 1.2296
\\
$2/1$
& 5.8168
& 2.669
& 1.2951
\\
\hline
\end{tabular}
\end{center}
\end{table}





\bibliographystyle{alpha}
\bibliography{squatJump_JBYBKM_2013_arXiv}

\end{document}